\documentclass{amsart}

\usepackage{amssymb}
\usepackage{mathtools}
\usepackage[svgnames]{xcolor}

\usepackage[unicode,
  colorlinks=true,
  linktocpage=true,
  citecolor=OliveDrab,
  linkcolor=DarkMagenta,
  linkcolor=DarkMagenta,
  pdftitle={The universal property of the decomposition space of quasisymmetric functions},
  pdfauthor={Philip Hackney and Joachim Kock and Jan Steinebrunner},
  pdfkeywords={simplicial object, quasisymmetric function, decomposition space, 2-Segal space, coalgebra}
  ]{hyperref}

\usepackage{enumitem}
\usepackage{tikz-cd}
\usepackage[capitalise,noabbrev]{cleveref}
\usepackage{microtype}
\usepackage[T1]{fontenc}
\usepackage{eucal}

\newcommand{\inert}{\mathbf{I}}
\newcommand{\contr}{\mathbf{C}}
\newcommand{\proto}{proto}

\newcommand{\yora}{a}
\newcommand{\xora}{a}

\newcommand{\gencopatt}{\mathbf{\Psi}}

\newcommand{\simpl}{\mathbf{\Delta}}
\newcommand{\fun}{\mathsf{Fun}}
\newcommand{\op}{\textup{op}}
\DeclareMathOperator{\id}{id}
\DeclareMathOperator*{\colim}{colim}
\DeclareMathOperator{\map}{map}
\DeclareMathOperator{\ob}{Ob}

\def\overarrow#1{{\vec{#1}}}
\def\nondeg{\overarrow}

\newcommand{\spaces}{\mathcal{S}}
\newcommand{\basecat}{\mathcal{E}}

\newcommand{\Span}{\mathsf{Span}}

\newcommand{\cat}{\mathsf{Cat}}
\newcommand{\catinf}{\mathsf{Cat}_\infty}
\newcommand{\scat}{\mathsf{sCat}}
\newcommand{\disc}{\operatorname{disc}}
\newcommand{\slicel}[2]{\vphantom{#2}^{{#1}/}\mkern-2mu{#2}}

\makeatletter
\newcommand*\wt[1]{\mathpalette\wthelper{#1}}
\newcommand*\wthelper[2]{%
        \hbox{\dimen@\accentfontxheight#1%
                \accentfontxheight#11.3\dimen@
                $\m@th#1\widetilde{#2}$%
                \accentfontxheight#1\dimen@
        }%
}

\newcommand*\accentfontxheight[1]{%
        \fontdimen5\ifx#1\displaystyle
                \textfont
        \else\ifx#1\textstyle
                \textfont
        \else\ifx#1\scriptstyle
                \scriptfont
        \else
                \scriptscriptfont
        \fi\fi\fi3
}
\makeatother

\newcommand{\Q}{\mathbb{Q}}
\newcommand{\Tsimpl}{\wt\simpl}
\newcommand{\Tord}{\Tsimpl}
\newcommand{\Tcontr}{\wt\contr}
\newcommand{\Tinert}{\wt\inert}
\newcommand{\gr}{\mathbf{R}}
\newcommand{\Tgr}{\wt\gr}
\newcommand{\dirgraph}{\gr}
\newcommand{\mdirgraph}{\Tgr}

\newcommand{\Tob}{\textup{\color{Indigo}T}}
\newcommand{\sigr}{{\color{Indigo}\varphi}}
\newcommand{\sigl}{{\color{Indigo}\varsigma}}

\newcommand{\lilh}{\mathsf{h}}
\newcommand{\lili}{\mathsf{i}}
\newcommand{\lilj}{\mathsf{j}}
\newcommand{\lilk}{\mathsf{k}}
\newcommand{\lilp}{\mathsf{p}}
\newcommand{\lilu}{\mathsf{u}}
\newcommand{\lilw}{\mathsf{w}}

\newcommand{\actrm}{\textup{act}}
\newcommand{\intrm}{\textup{in}}

\newcommand{\pre}{\mathcal{P}}
\newcommand{\culf}{\textup{culf}}
\newcommand{\ikeo}{\textup{ikeo}}
\newcommand{\cons}{\textup{cons}}
\newcommand{\std}{\textup{std}}

\newcommand{\QSym}{\operatorname{QSym}}
\newcommand{\N}{\mathbb{N}}
\newcommand{\tensor}{\otimes}

\newcommand{\dcons}[1]{\mathsf{Cons}(#1)}
\newcommand{\dgr}[1]{\mathsf{Gr}(#1)}

\usepackage[old]{old-arrows}
\usepackage[only,mapsfromchar,mapsfrom]{stmaryrd}
\newcommand{\ract}{\varrightarrow\mapsfromchar}
\newcommand{\rint}{\rightarrowtail}
\tikzset{
  act /.tip = >|
}

\tikzcdset{
  ds/.style={execute at begin picture={\everymath\expandafter{\the\everymath\displaystyle}}}
}

\newcommand{\drpullback}{\arrow[phantom]{dr}[very near start,description]{\lrcorner}}

\newcommand{\name}[1]{\ulcorner \!#1\! \urcorner}

\newtheorem{theorem}{Theorem}[section]
\newtheorem{letterthm}{Theorem}
\newtheorem*{untheorem}{Theorem} 
\newtheorem{lemma}[theorem]{Lemma}
\newtheorem{proposition}[theorem]{Proposition}
\newtheorem{corollary}[theorem]{Corollary}
\theoremstyle{definition}
\newtheorem{definition}[theorem]{Definition}
\newtheorem{example}[theorem]{Example}
\newtheorem{notation}[theorem]{Notation}

\newtheorem{convention}[theorem]{Convention}
\newtheorem{remark}[theorem]{Remark}

\newtheoremstyle{noname}
  {8pt plus 2pt minus 4pt} 
  {8pt plus 2pt minus 4pt} 
  {}                       
  {}                       
  {\bfseries}              
  {.}                      
  {5pt plus 1pt minus 1pt} 
  {\thmnumber{#2}\thmnote{.\ #3}} 

\theoremstyle{noname}
\newtheorem{blanko}[theorem]{}

\makeatletter
\g@addto@macro\th@plain{\let\thmhead\thmhead@plain}
\g@addto@macro\th@definition{\let\thmhead\thmhead@plain}
\g@addto@macro\th@remark{\let\thmhead\thmhead@plain}
\makeatother

\AddToHook{env/lemma/begin}{\crefalias{theorem}{lemma}}
\AddToHook{env/proposition/begin}{\crefalias{theorem}{proposition}}
\AddToHook{env/corollary/begin}{\crefalias{theorem}{corollary}}
\AddToHook{env/definition/begin}{\crefalias{theorem}{definition}}
\AddToHook{env/example/begin}{\crefalias{theorem}{example}}
\AddToHook{env/notation/begin}{\crefalias{theorem}{notation}}
\AddToHook{env/construction/begin}{\crefalias{theorem}{construction}}
\AddToHook{env/convention/begin}{\crefalias{theorem}{convention}}
\AddToHook{env/remark/begin}{\crefalias{theorem}{remark}}
\apptocmd{\appendix}{\crefalias{section}{appendix}}{}{}

\AtBeginDocument{%
   \def\MR#1{}
}

\begin{document}

\title[Universal property of quasisymmetric functions]{The universal property of the decomposition space of quasisymmetric functions}

\author{Philip Hackney}
\address{Department of Mathematics, University of Louisiana at Lafayette, USA}
\email{philip@phck.net}

\author{Joachim Kock}
\address{Departament de Matem\`atiques, Universitat Aut\`onoma de Barcelona, Spain and Centre de Recerca Matem\`atica (Barcelona), Spain}
\email{joachim.kock@uab.cat}

\author{Jan Steinebrunner}
\address{University of Oxford, UK}
\email{jan.steinebrunner@maths.ox.ac.uk}

\thanks{ This work was supported by a grant from the Simons Foundation
(\#850849, PH). The first author was partially supported by the Louisiana
Board of Regents through the Board of Regents Support fund
LEQSF(2024-27)-RD-A-31. The second and third authors were fully funded by 
the Independent Research Fund Denmark under grant number 
10.46540/3103-00099B,
and acknowledge also support from the Danish National Research
Foundation through the Copenhagen Centre for Geometry and Topology (DNRF151).
The second author acknowledges support from research grants PID2020-116481GB-I00
and PID2024-158573-NBI00
(AEI/FEDER, UE) of Spain and grant 2025-SGR-01322 of Catalonia, and through the
Severo Ochoa and Mar\'ia de Maeztu Program for Centers and Units of Excellence in
R\&D grant number CEX2020-001084-M }

\begin{abstract}
Aguiar, Bergeron, and Sottile showed that the coalgebra of quasi\-symmetric functions is terminal among connected graded coalgebras with a zeta function.
After developing the notions of gradings and zeta functions for decomposition spaces, we establish an objective homotopic version of this result. 
Namely, the decomposition space $Q$ of quasisymmetric functions is terminal: for any conservatively graded decomposition space $X$ with zeta function $F$, there is a unique graded ikeo-culf span $X \leftarrow W \rightarrow Q$ inducing $F$ from the special zeta function on $Q$.
More generally, for gradings in
any complete decomposition space $M$ with indecomposable identities, the terminal object is $J(M)$, the free decomposition space on the nondegenerate simplices of $M$; the quasisymmetric case is $Q=J(B\mathbb{N})$.
\end{abstract}

\maketitle

\section{Introduction}

\begin{blanko}[Quasisymmetric functions]
  Quasisymmetric functions were introduced and studied by
Gessel~\cite{Gessel:84} as generating functions for $P$-partitions,
based on earlier work of Stanley~\cite{Stanley:1972} where the quasisymmetric
functions already appear implicitly. That $\QSym$ is in fact a Hopf algebra was
shown independently by Malvenuto--Reutenauer~\cite{Malvenuto-Reutenauer} and
Ehrenborg~\cite{Ehrenborg}.\footnote{Gessel~\cite{Gessel:84} described the so-called
inner comultiplication law (for which the power-sum symmetric functions are
group-like). The comultiplication relevant in the present work is the outer
comultiplication, first described by Malvenuto and
Reutenauer~\cite{Malvenuto-Reutenauer} (for which the power-sum symmetric
functions are primitive). This comultiplication (the only one we consider) is
described below (deconcatenation).} We refer to Hazewinkel~\cite{Hazewinkel:2003},
Grinberg--Reiner~\cite{Grinberg-Reiner:1409.8356} or 
Luoto--Mykytiuk--van Willigenburg~\cite{Luoto-Mykytiuk-vanWilligenburg}
for background.

Since then, quasisymmetric functions have turned up throughout algebraic combinatorics, as well
as in applications areas, for example in representation theory
\cite{Hoffman:9907173}, \cite{Krob-Thibon:NSym-IV},
\cite{Billera-Brenti:0710.3965}, in graph theory
\cite{Shareshian-Wachs:1405.4629}, \cite{Alexandersson-Sulzgruber:1807.02460},
\cite{Mason:1810.07148}, in topology~\cite{Baker-Richter:0605743}, in
stochastic analysis \cite{Diehl-EbrahimiFard-Tapia:1906.05823}, and free
probability~\cite{Mastnak-Nica:TAMS2010}, just to mention a few.

The quickly increasing number of occurrences of quasisymmetric generating
functions associated to combinatorial structures, led to rather systematic
constructions such as those given by Bergeron, Mykytiuk, Sottile, and van
Willigenburg~\cite{Bergeron-Mykytiuk-Sottile-vanWilligenburg}, and culminated with
the
following theorem of Aguiar, Bergeron, and 
Sottile~\cite{Aguiar-Bergeron-Sottile}:

\begin{untheorem}
The Hopf algebra (or just coalgebra) of quasisymmetric functions is the terminal
object in the category of connected graded Hopf algebras (or just coalgebras) 
with a zeta function.
\end{untheorem}

It explains quite clearly why quasisymmetric functions turn up in so many contexts.
While in retrospect the result is not surprising, the formula for this universal map
is rather involved (see~\ref{ABS} further on), and it encodes quite a lot
of information.
\end{blanko}

In this paper we establish an objective homotopic version of this result, 
stating the corresponding $\infty$-categorical  universal property of the decomposition space
of quasi\-symmetric functions:
\begin{letterthm}\label{intro theorem A}
The decomposition space $Q$ of quasisymmetric functions, with its special
zeta function, is the terminal object in the category whose objects are
conservatively graded decomposition spaces with zeta functions, and whose
morphisms are graded ikeo-culf spans compatible with zeta functions.
\end{letterthm}
We note
that multiplicative structure plays no role in this theorem, which only concerns the
coalgebraic setting; we hope to address the Hopf algebra case in future work.

\Cref{intro theorem A} is part of a general program to develop algebraic combinatorics over
decomposition spaces, and an important secondary goal of the paper is to set up the
general theory required.
Every phrase in this theorem needs to be explained (including \emph{grading} and
\emph{zeta functions}), and we do so momentarily, but first we give a reminder of the
way that the theory of decomposition spaces interacts with combinatorics. 

\begin{blanko}[Decomposition spaces -- objective algebraic combinatorics]
  Algebraic combinatorics has led to the solution of many difficult
  combinatorial problems, and has been the vehicle for deep interactions with
  geometry, probability theory, representation theory,  quantum algebra, 
  mathematical physics, and so on. At the same time, bijective
  proofs are well appreciated in combinatorics as an expression of deeper understanding
  than algebraic identities alone can give. It is thus
  natural to try to objectify algebraic combinatorics, which means to mimic its
  methods directly at the level of combinatorial objects rather than with the algebra
  they generate, so that algebraic identities arise by taking
  cardinality of bijections of sets (or homotopy cardinality of equivalences of
  groupoids or spaces). A blueprint for such objectification is
  Joyal's theory of species~\cite{Joyal:1981} which showed how the basic 
  constructions and arguments with
  generating functions can be performed directly at the level of combinatorial
  species.
  
  Lawvere and Menni~\cite{Lawvere-Menni} pioneered the program of doing combinatorial
  Hopf and coalgebras at the objective level, based on Leroux's theory of M\"obius
  categories~\cite{Leroux:1976}.
  Decomposition spaces \cite{GKT1,GKT2,GKT3} (a.k.a.\ $2$-Segal
  spaces \cite{DyckerhoffKapranov:HSS}) are a wide-ranging generalization of these
  ideas: they are certain simplicial spaces $X\colon \simpl^\op\to\spaces$ which generalize (nerves of) 
  posets, monoids,
  categories, and operads, and they allow the construction of incidence (co)algebras
  and a M\"obius inversion principle. Like in all these cases it is $X_1$ (simplicial 
  degree $1$) that generates the underlying vector space of the coalgebra, while the 
  remaining simplicial degrees serve to express the structure.
  The objective analog of the category of vector spaces is the category of
  spans (of $\pi$-finite spaces)~\cite{Galvez-Kock-Tonks:1602.05082}, which can thus be seen as an 
  intermediate level between that of decomposition spaces and that of vector spaces; the classical incidence (co)algebras then arise 
  from
  a cardinality functor from spans to vector spaces.
  The relationship between the decomposition-space level and the level of coalgebras in spans is the topic of some deep theorems: First
  it was shown by Stern~\cite{Stern:1905.06671} that
  every (co)algebra in spans is the incidence coalgebra of a decomposition space, and 
  it was recently shown by G\"odicke~\cite{Godicke:IC2SS}
  that every (co)algebra homomorphism in spans 
  comes from an ikeo-culf span of decomposition spaces (detailed in \ref{culf+ikeo} below). If a 
  co-, bi- or Hopf algebra is defined by combinatorial data, it can usually
  be realized as a co-, bi- or Hopf algebra in spans, and therefore has a
  decomposition-space model which we regard as its objective counterpart.
  Throughout this paper we work at the level of decomposition spaces.
      
  While many basic operations with   combinatorial
  (co)algebras can be performed at this level of
  generality~\cite{Galvez-Kock-Tonks:1612.09225}, it is not reasonable to expect
  that all of algebraic combinatorics can be recast in this 
  language -- in the big picture perhaps only a small fraction,
  but it is
  still very fruitful to see how far combinatorics can be stretched into 
  algebra, and enlightening to see constructions and results that are usually
  very algebraic in appearance reformulated in purely combinatorial (and 
  categorical) terms. A surprising number of conditions and arguments end up
  formulated in terms of pullbacks.

  The theory of decomposition spaces is staged in simplicial homotopy theory and has a
  strong categorical flavor. Homotopy theory and $\infty$-categories come in for two
  reasons: one is that the theory is simplicial and categorical from the outset, so it
  naturally drifts into homotopy theory. But the combinatorial reason why homotopy
  theory comes in is that combinatorial objects have symmetries, and to take them into
  account natively requires groupoids: symmetries are already homotopies.
    
\end{blanko}

\begin{blanko}[Gradings and zeta functions]
One contribution of this paper is to develop gradings and zeta functions for decomposition spaces.
  At the level of decomposition spaces, a grading is a simplicial map to
  $B\N$, the nerve of the monoid $(\N,+)$. In the Aguiar--Bergeron--Sottile
  theorem, the graded coalgebras are \emph{connected} but we use the weaker
  \emph{proto-connected}; for decomposition spaces this means the grading is 
  conservative (\cref{proto-connected = conservative}).

  A {\em zeta function} on a combinatorial coalgebra is a linear
  functional that agrees with the counit in degree $0$. 
  At the level of decomposition spaces, a zeta function is
  given by a single map $X_1 \leftarrow F$ (we write it left-pointing 
  because it is really the left-hand part of a span $X_1 \leftarrow F \to 
  1$), and the condition on
  degree $0$ translates into the requirement that we have a
  pullback diagram (cf.\ \ref{zeta decomp})
  \begin{equation*}
  \begin{tikzcd}
  X_0 \ar[d, equal] \drpullback \ar[r] & F \ar[d]  \\
  X_0 \ar[r, "s_0"'] & X_1 .
  \end{tikzcd}
  \end{equation*}

\end{blanko}

\begin{blanko}[The decomposition space $Q$]
  The coalgebra $\QSym$ (like most combinatorial coalgebras) is the
  incidence coalgebra of a decomposition space~\cite{GKT1}, namely a
  certain decomposition space $Q$ which in simplicial degree $1$ has 
  $Q_1 = (\N_+)^*$, the
  set of words in positive integers~\cite{GKT:QSym},
  \cite{Hackney-Kock:2210.11192}, cf.~\ref{decomp qsym} below. 

The relevant zeta function on $\QSym$ is not the \emph{standard}
one (which takes every basis element to $1$), but a \emph{special} one,
taking value $1$ on all words of length $0$ and $1$ and value $0$ on longer
words. At the level of the decomposition space $Q$, this special zeta
function is given by the map $(\N_+)^*\leftarrow \N$ that takes positive
number to the corresponding one-letter word and $0$ to the empty word. (In
\ref{zeta of free} we propose a general notion of ``special'' zeta function.)

\end{blanko}

We have now discussed the main ingredients in \cref{intro theorem A}. It turns out the 
proof does not depend much on $B\N$, but works the same 
if $B\N$ is replaced by a more general decomposition space $M$ as 
recipient of the grading:

\begin{letterthm}\label{intro theorem B}
  For any complete decomposition space $M$ with indecomposable identities,
  and for any conservatively $M$-graded decomposition space $X$ with a zeta
  function $X_1 \leftarrow F$, there is a unique $M$-graded ikeo-culf span $X
  \leftarrow W \to J(M)$ that induces $F$ from the special zeta function on
  $J(M)$.
\end{letterthm}

Here the terminal object $J(M)$ is the free decomposition space of nondegenerate simplices in $M$,
which in the case $M=B\N$ is the usual decomposition space of quasisymmetric
functions: $Q = J(B\N)$. 
  This generalization reveals that the essential feature of $Q$, to which it
  owes its universal property, is the fact that it is the free decomposition
  space on the nondegenerate part of the recipient of the notion of grading.

\begin{blanko}[The $J$-construction] 
  The intuitive description of the $J$-construction is a
  two-step process: first replace $M$ with its collection of
  nondegenerate simplices $\nondeg M$,
  which is a presheaf over the subcategory $\inert \subset \simpl$ of 
  inert maps
  (i.e.\ interval inclusions), and then left Kan extend along
  the inclusion $\lilj$ to obtain the simplicial space 
  $J(M) =\lilj_! \nondeg M$. 
  However, the first step is not functorial in all maps, and works only for complete
  and stiff simplicial spaces.
  Nevertheless, the composite, $J$, is
  actually functorial! To see this, we show that $J$ can be given another
  description, valid for general simplicial spaces. This exploits an
  intermediate category
$$
\inert \longrightarrow \contr \stackrel \lilk \longrightarrow \simpl,
$$
the category of contractive maps
(generated by degeneracies and outer face maps). 
We show (\cref{cor JX identification}) that 
$$
J(M) = \lilj_! \nondeg M = \lilk_! \lilk ^* M,
$$
and clearly $J\coloneq \lilk_!\lilk^*$ works without any assumptions on the simplicial 
space at all, and is functorial in all simplicial maps.
The indexing category $\contr$ plays an important role, as a natural 
setting for the notions of complete and stiff (cf.~\cite[2.6 and 
\S4]{GKT2}).
\end{blanko}

\begin{blanko}[Overview]
  The following brief section, \cref{sec universal coalgebra qsym}, gives a touch more
  detail on quasisymmetric functions and the original Aguiar--Bergeron--Sottile
  theorem. It is followed by \cref{sec decomp ikeo culf q}, which gives background on
  the decomposition space analog of the coalgebra of quasisymmetric functions via
  free decomposition spaces, and on G\"odicke's theorem, which gives the correct
  notion of map of decomposition spaces for modeling maps of coalgebras. Sections
  \ref{sec:grading} and \ref{sec:zeta} are new, and develop the appropriate notions of
  grading and zeta functions in this setting. 
  These four sections expand on
  central notions in \cref{intro theorem A} and \cref{intro theorem B}.
  The final section, \cref{sec finiteness}, imposes additional finiteness conditions
  and interprets \cref{intro theorem A} in terms of the Aguiar--Bergeron--Sottile
  formula for the unique map to $\QSym$. All of these sections take place over the
  $\infty$-category of spaces, $\spaces$, and can be seen as grounding the theorems in
  concrete combinatorics.

The intermediate sections are more technical, and can largely be read independently of sections \ref{sec universal coalgebra qsym}--\ref{sec:zeta} and \ref{sec finiteness}.
These take place over a somewhat general base category or $\infty$-category $\basecat$ satisfying the hypotheses in \ref{base category}, and culminate in a proof of \cref{intro theorem B} (as \cref{main theorem}).
This is the natural generality where our proof of the main theorem holds, but the reader is welcome to replace $\basecat$ by the 1-category of sets or the $\infty$-category of spaces throughout.
\Cref{sec shape categories} introduces the indexing category $\Tsimpl$, which encodes zeta functions, and two important subcategories $\Tinert \subset \Tcontr \subset \Tsimpl$, and explains factorization systems on these. 
It is followed by \cref{sec object prop} which explains a number of relevant conditions on presheaves over these categories,
and \cref{sec map properties} which explains several important kinds of maps (and the definition of the $\infty$-category $\dgr{M}$ of decomposition spaces conservatively graded in $M$).

The next five sections concern the primary ingredients used in the proof of the main theorem.
\Cref{section standard} concerns \emph{standard} presheaves, namely those whose zeta function $X_1 \leftarrow F$ is an equivalence, describes the universal standard approximation to a general $\Tcontr$-presheaf, and establishes a key contractibility lemma used in the proof of the main theorem.
We next give explicit formulas for left Kan extension along the inclusions $\Tinert \subset \Tcontr$ and $\Tinert \subset \Tsimpl$ in \cref{sec lke along jkl}.
Interaction between these left Kan extensions and culf/conservative maps is the subject of \cref{sec culf conserv lke}, which culminates in \cref{thm restricted slices} giving several equivalences between restricted slice categories.
We then characterize which $\Tcontr$-presheaves $A$ arise as left Kan extensions along $\Tinert \subset \Tcontr$ in \cref{sec nondeg elts}, and 
in particular explain that there is at most one $\nondeg A \in \pre(\Tinert)$ whose left Kan extension is $A$; an alternate title for this section is ``Homotopical Eilenberg--Zilber Theory''.
\Cref{sec:ikeo} explains the tight relationship between left Kan extension along $\Tcontr \subset \Tsimpl$ and ikeo maps.
Finally, \cref{sec main thm proof} assembles the proof of the main theorem.
\end{blanko}

\begin{blanko}[Base category assumptions]\label{base category}
In sections \ref{sec object prop}--\ref{sec main thm proof} we will work over a sufficiently nice base category or $\infty$-category $\basecat$ under the hypotheses below.
Suitable bases include any $\infty$-topos (e.g.\ the $\infty$-category of spaces $\spaces$) or any ordinary Grothendieck 1-topos (e.g.\ the category of sets).
\begin{itemize}[left=0pt]
\item \textbf{Limits and colimits:} We assume that $\basecat$ is closed under small limits and colimits.
\item \textbf{Colimits are universal:} For each map $X \to Y$ in $\basecat$, the base change map $\basecat_{/Y} \to \basecat_{/X}$ preserves colimits.
\item \textbf{Extensivity:} For any small family $(X_a)$ of objects of $\basecat$, the canonical map
\[ \begin{tikzcd}
\basecat_{/\sum_a X_a} \rar{\simeq} & \displaystyle \prod_a \basecat_{/X_a} 
\end{tikzcd} \]
is an equivalence.
\end{itemize}
In these sections, including in the main theorem, $\pre(\mathcal{C})$ will denote $\fun(\mathcal{C}^\op, \basecat)$.
In sections \ref{sec universal coalgebra qsym}--\ref{sec:zeta} and \ref{sec finiteness} we specialize to the case $\basecat = \spaces$.
\end{blanko}

\subsection*{AI Statement:} 
The main results of this paper were obtained in 2024 and presented by the second author at CT2024. 
Late in the revision of a complete manuscript, while thinking about \cref{sec:ikeo} and how to weaken its assumption on the base category $\basecat$ (originally assumed to be an $\infty$-topos), we consulted a large language model about this question. 
It suggested a different approach, which the authors then used to redevelop and simplify \cref{sec:ikeo}. 
Apart from that, no AI was used.

\section{The universal property of the coalgebra of quasisymmetric functions}\label{sec universal coalgebra qsym}

\begin{blanko}[Quasisymmetric functions, classically]
  Recall (for example from \cite{Stanley:volII}) that quasisymmetric functions
  are certain power series of bounded degree in countably many variables
  $x_1,x_2,\ldots$: for each word $w$ of length $k$ in the alphabet $\N_+$ there
  is a quasi-symmetric function
  $$
  M_w \coloneq \sum_{i_1 < i_2 <\cdots < i_k} x_{i_1}^{w_1} x_{i_2}^{w_2} \cdots
  x_{i_k}^{w_k} ,
  $$
  and the ring of quasisymmetric functions is the linear span of these
  so-called {\em monomial} quasisymmetric functions $M_w$. The set $Q_1 =
  (\N_+)^*$ of all words in $\N_+$ thus indexes this {\em monomial basis}. For
  example, to a word such as $(23114)$ corresponds the quasisymmetric function
  $M_{23114} = \sum_{i_1 < i_2 <\cdots < i_5} x_{i_1}^2 x_{i_2}^3 x_{i_3}^1
  x_{i_4}^1 x_{i_5}^4$. Although the actual power series are important in many
  applications, the algebraic structure is governed by the combinatorics of
  words in $\N_+$ (in the $M$-basis and in several related bases), and in the
  following we shall not need the actual power series. In particular the
  coalgebra structure of $\QSym$ can be described purely in terms of words: the
  comultiplication is simply word splitting (deconcatenation).
\end{blanko}

\begin{blanko}[Gradings and the special zeta function]
  \label{Q-gradings-and-special-zeta}
  There are two important gradings on $\QSym$: the length grading (sending
  $M_{(n_1,\ldots,n_k)}$ to $k$) and the total grading (sending
  $M_{(n_1,\ldots,n_k)}$ to $\sum n_i$). The coalgebra $\QSym$ is connected for
  both gradings. It is with respect to the total grading that $\QSym$ has the
  universal property.
  
  The {\em special zeta function} on $\QSym$ is the linear functional 
  $\zeta_{\operatorname{special}} \colon \QSym \to \Q$
  that takes value $1$ on all words of length $0$
  and $1$, and takes value $0$ on all longer words
  (cf.~\cite{Aguiar-Bergeron-Sottile}).
\end{blanko}

\begin{blanko}[Aguiar--Bergeron--Sottile formula]\label{ABS}
  If $(C,\phi)$ is a graded coalgebra with zeta function, then the universal coalgebra
  map $(C,\phi)\to(\QSym, 
  \zeta_{\operatorname{special}})$ is given by
  $$
  x \ \mapsto \ \sum_{r\geq 0} \ \sum_{x_{(1)}\otimes\dots\otimes x_{(r)}\in \Delta^{(r)}(x)}
  (\phi-\epsilon)^{\otimes r}(x_{(1)}\otimes\dots\otimes x_{(r)}) \cdot 
  M_{\deg(x_{(1)})\cdots \deg(x_{(r)})} .
  $$
  The sums together are thus over all possible comultiplications of $x$ (including 
  $\Delta^{(0)} =\epsilon$ and $\Delta^{(1)} = \id$), and 
  the corresponding monomial is indexed by the degree vector of the factors of the 
  comultiplication. The coefficient is the value of the $r$-fold tensor product of 
  $\phi-\epsilon$. If $x$ is of degree $k$, then all comultiplications longer than $k$
  will contain a tensor factor of degree $0$, and then $\phi-\epsilon$ will vanish on 
  this factor, killing the whole term. In this way, for fixed $x$ the sum is finite,
  and at the same time the vanishing of $\phi-\epsilon$ on degree $0$ means that 
  the number $0$ never appears in the degree vector, which is therefore a valid index for the 
  monomial.
\end{blanko}

To see how this formula induces $\phi$ from 
$\zeta_{\operatorname{special}}$,
apply $\zeta_{\operatorname{special}}$ to the whole formula, and check that the 
result is $\phi(x)$. Only $r=1$ and $r=0$ 
can give any contribution, since $\zeta_{\operatorname{special}}$ vanishes on 
higher-degree generators. So for the purpose of this calculation the relevant part of 
the right-hand side is
$$
  \epsilon(x) \cdot 
  M_{\{\}} 
  +
  (\phi-\epsilon)(x) \cdot 
  M_{\deg(x)} .
$$
Since $\zeta_{\operatorname{special}}$ takes value $1$ on these generators,
we are left with  
$$\epsilon(x)
  +
  (\phi-\epsilon)(x) = \phi(x)
$$
as required.

The $r>1$ terms did not enter this calculation, but are required to
make the map a graded coalgebra homomorphism. 
The proof of the universal property given by Aguiar, 
Bergeron and Sottile uses duality, and argues in terms of noncommutative 
symmetric functions (the graded dual of $\QSym$).

The proof we give works with more global objects, and in particular produces
the span
whose expansion in coordinates is the above formula. We will not see the formula
again until at the very end of the paper, for comparison.

\section{Decomposition spaces, ikeo-culf spans, and quasisymmetric functions}\label{sec decomp ikeo culf q}

This section recalls some background and notation.

\begin{blanko}[Decomposition spaces and incidence coalgebras]
  The simplex category $\simpl$ has an active-inert factorization system.
  The active maps, forming the subcategory $\simpl^\actrm \subset \simpl$, are 
  generated by the degeneracy maps and the inner face maps. The inert maps,
  forming the subcategory $\inert \coloneq \simpl^\intrm$, are spanned by the outer face maps.
  A {\em decomposition space}~\cite{GKT1} (also called {\em $2$-Segal 
  space}~\cite{DyckerhoffKapranov:HSS})
  is a presheaf $X\in \pre(\simpl)$ that sends
  active-inert pushouts to pullbacks.

  The {\em incidence coalgebra} of a decomposition space $X$ is the slice category
  $\spaces_{/X_1}$ whose comultiplication is given by the span
\[ \begin{tikzcd}
X_1 & \lar[swap]{d_1} X_2 \rar{(d_2,d_0)} & X_1 \times X_1
\end{tikzcd} \]
  and whose counit is given by the span
  $$
  X_1 \stackrel{s_0}\leftarrow X_0 \to 1.
  $$
  The decomposition-space axioms ensure that this coalgebra is 
  coassociative~\cite{GKT1}.
\end{blanko}

\begin{blanko}[Ikeo-culf spans and coalgebra homomorphisms]\label{culf+ikeo}
There are two important notions of maps between decomposition spaces: the culf maps and the ikeo maps (\cref{sec map properties} or \cite{GKT1}). 
Culf maps induce coalgebra homomorphisms
  covariantly (\cite[Lemma~8.2]{GKT1}), while ikeo maps induce coalgebra homomorphisms contravariantly (\cite[Proposition~8.6]{GKT1}).

  Altogether, an ikeo-culf span
\[ \begin{tikzcd}
Y  & \lar["\text{ikeo}"']  W \rar["\text{culf}"] &  X
\end{tikzcd} \]
  induces a coalgebra homomorphism from the incidence coalgebra of $Y$ to the 
  incidence coalgebra of $X$.
It had been expected for a long time
  (first formulated by Joyal) that {\em all} coalgebra homomorphisms arise in this way. 
This was recently proved by G\"odicke:
\end{blanko}

\begin{theorem}
  [{G\"odicke~\cite[Theorem 1.1]{Godicke:IC2SS}}]
  Let $\basecat$ be an $\infty$-category with finite limits.
  The incidence-coalgebra construction constitutes an equivalence
  between the $\infty$-category of decomposition spaces in $\basecat$ with ikeo-culf spans,
  and the $\infty$-category of coalgebras in $\Span(\basecat)$.
\end{theorem}
In view of this result, the uniqueness in our main theorem (a statement about 
ikeo-culf spans of decomposition spaces)
can also be interpreted as a uniqueness in the category of coalgebras in spans.

\begin{blanko}[Free decomposition spaces]
  Given a $\inert$-presheaf $A$, that is,
  a diagram $A_0 \leftleftarrows A_1 \leftleftarrows A_2 \leftleftarrows \cdots$, 
  left Kan extension along the inclusion $\lilj\colon \inert\hookrightarrow \simpl$ defines a simplicial space
  $\lilj_!(A) \in \pre(\simpl)$ which is automatically a decomposition 
  space~\cite[Corollary 2.3.3]{Hackney-Kock:2210.11192}. It is described explicitly as
  \begin{eqnarray*}
    \lilj_!(A)_0 & = & A_0 \\
    \lilj_!(A)_1 & = & \sum_{n\in \N} A_n
  \end{eqnarray*}
  and generally
  \begin{eqnarray*}
    \lilj_!(A)_k &=& \!\!\!\sum_{\alpha\colon [k]\ract [n]} A_n  .
  \end{eqnarray*}
  Decomposition spaces that arise in this way are called {\em free}.
\end{blanko}

\begin{blanko}[Free decomposition spaces of nondegenerate simplices]\label{J}
  A decomposition space is {\em complete} \cite{GKT2} when the degeneracy maps are monomorphisms.
  In particular $s_0\colon X_0 \to 
  X_1$ then has a complement denoted 
  $\nondeg X_1$ consisting of {\em nondegenerate} $1$-simplices.
  Similarly one 
  defines $\nondeg X_n \subset X_n$ to be the complement of all the 
  degeneracy maps, which is the subspace consisting of the 
  $n$-simplices all of whose principal edges are nondegenerate $1$-simplices (\cite[Corollary 2.16]{GKT2}).

  One may then form the {\em free decomposition space $J(X)$ of nondegenerate simplices}
  \cite{GKT:QSym}, by first constructing the $\inert$-presheaf 
  $\nondeg X_0 \leftleftarrows \nondeg X_1 \leftleftarrows \nondeg X_2 
  \leftleftarrows \cdots$ 
  and then setting $J(X) \coloneq \lilj_!(\nondeg X)$. 
  We have
  \begin{eqnarray*}
  J(X)_0 & = & X_0 \\
  J(X)_1 & = & \sum_{n\in \N} \nondeg X_n ,
  \end{eqnarray*}
  the set of all nondegenerate simplices in all dimensions.

  In higher simplicial dimension, it has {\em subdivided} nondegenerate 
  simplices. Precisely, a $k$-simplex of $J(X)$ is a diagram
  \[ \begin{tikzcd}
  \Delta^k  \rar[-act] &  \Delta^n \rar{\text{nondeg}} &[+0.2cm] X ,
  \end{tikzcd} \]
  where the second map determines the nondegenerate $n$-simplex in $X$ and the 
  active map $\Delta^k \ract \Delta^n$ gives the subdivision, whereby 
  the $n$-simplex is subdivided into $k$ ``stages''. 
   
  There is a canonical simplicial map
  $$
  X \leftarrow J(X)
  $$
  given in dimension $1$ by sending a nondegenerate simplex to its long edge. In
  dimension $k$ it sends a $k$-subdivided $n$-simplex to its $k$-simplex of
  stages. This map is always ikeo 
  (as we shall see in a more general
  viewpoint in \cref{sec:ikeo} (specifically \cref{lem ikeo counit})).
\end{blanko}

\begin{blanko}[The decomposition space of quasisymmetric functions]\label{decomp qsym}
  In the special case where the decomposition space $X$ is $B\N$, the nerve of the monoid
  $(\N,+)$, then the $\inert$-presheaf of nondegenerate simplices
  has $(\N_+)^n$ in degree $n$, and the free decomposition space
  $$
  J(B\N) = \lilj_!(\nondeg{B\N}) \eqcolon Q
  $$ 
  is the decomposition space of words in $\N_+$,
  whose
  incidence coalgebra is the coalgebra of quasisymmetric functions
  $\QSym$ in the $M$-basis~\cite{GKT:QSym}; see also~\cite[5.1.4]{Hackney-Kock:2210.11192}.
\end{blanko}

\section{Grading}

\label{sec:grading}

The importance of gradings of combinatorial coalgebras is well established.
Aguiar--Bergeron--Sottile~\cite{Aguiar-Bergeron-Sottile} even take
``combinatorial coalgebra'' to mean ``coalgebra equipped with a grading''. The
grading serves to perform inductive arguments, expressing how combinatorial
objects can be broken into smaller ones (i.e.~of lower degree). It is very
common to work furthermore with {\em connected} graded coalgebras (and this is
the setting of \cite{Aguiar-Bergeron-Sottile}), meaning that the zeroth graded 
piece is of
dimension $1$. The connectedness assumption gives a convenient basis for
the induction. An important expression of this is the case of bialgebras, where
connectedness ensures the existence of an antipode: a connected graded bialgebra
is a Hopf algebra. 

However, for most practical purposes it is just as good to demand the coalgebras to be
{\em \proto-connected}, meaning that degree zero is spanned by group-like elements --
this suffices to provide a basis for
inductions.\footnote{``Proto-connected'' was called ``not-quite-connected'' in
\cite{Kock:1411.3098}, \cite{Kock:1809.00941}, where it was argued that it is a
suitable setting for bialgebraic renormalization. Even in the case of bialgebras,
there is actually a useful notion of antipode also in the \proto-connected
case~\cite{Carlier-Kock:1807.11858}, although it does not have all the properties of
the full-blown antipode of the connected case.} 
From a \proto-connected graded coalgebra $C$, one gets a connected
graded coalgebra $C/I$ by quotienting by the graded coideal $I$ spanned by all $h-g$
for $h,g$ group-like. In this way ``connected'' and ``\proto-connected'' are not very
different, but \proto-connected coalgebras arise more frequently in nature.

In this paper we work with \proto-connected graded coalgebras, and we work with them 
in the guise of decomposition spaces, defined in \ref{proto} below. 
The main result of this section is that a proto-connected grading on a decomposition space $X$ is precisely a conservative simplicial map $X\to B\N$ (\cref{proto-connected = conservative}).
The impatient reader could take this as a definition and skip the 
rest of this
section, but it is important for us to establish that it is correct notion from
the viewpoint of algebraic combinatorics. It is also striking that yet another 
important notion from combinatorics ends up as a pullback condition.

Before coming to grading, we need the notion of group-like element. Recall
(e.g.~from \cite{Sweedler}) that an element $x$ in a coalgebra is group-like
when $\Delta(x)=x\tensor x$ and $\epsilon(x)=1$. At the level of 
decomposition spaces $X$, ``elements'' are  objects in the slice $\spaces_{/X_1}$,
that is maps $A\to X_1$.
This motivates the following
definition.

\begin{blanko}[Group-like elements]
  A generalized $1$-simplex $A \to X_1$ in a decomposition space $X$ is called {\em group-like}
  when the pullback along an active map $X_n \to X_1$ produces the identity on $A$ 
  like this:
  \[
  \begin{tikzcd}
  A \drpullback \ar[d, equal] \ar[r] & X_n \ar[d, -act]  \\
  A \ar[r] & X_1 .
  \end{tikzcd}
  \]
  In view of the decomposition-space axioms, it is enough to check the two basic 
  active maps
  \[
  \begin{tikzcd}
  A \drpullback \ar[d, equal] \ar[r] & X_0 \ar[d, -act, "s_0"]  \\
  A \ar[r] & X_1
  \end{tikzcd}
  \qquad
  \begin{tikzcd}
  A \drpullback \ar[d, equal] \ar[r] & X_2 \ar[d, -act, "d_1"]  \\
  A \ar[r] & X_1  .
  \end{tikzcd}
  \]
In particular, group-like implies degenerate.
\end{blanko}

\begin{remark}
  To see what this has to do with group-like elements in the classical sense, consider the case $A=1 \xrightarrow{a} X_1$. Then the square on the left expresses $\epsilon(\xora)=1$, while
  at first sight it looks like the square on the right only ensures that $\Delta(\xora) = 
  \xora'\tensor \xora''$ for possibly different $\xora'$ and $\xora''$. But in fact we do have
  $\Delta(\xora)=\xora\tensor \xora$, as can be seen from the diagram (with $i=0,1$)
  
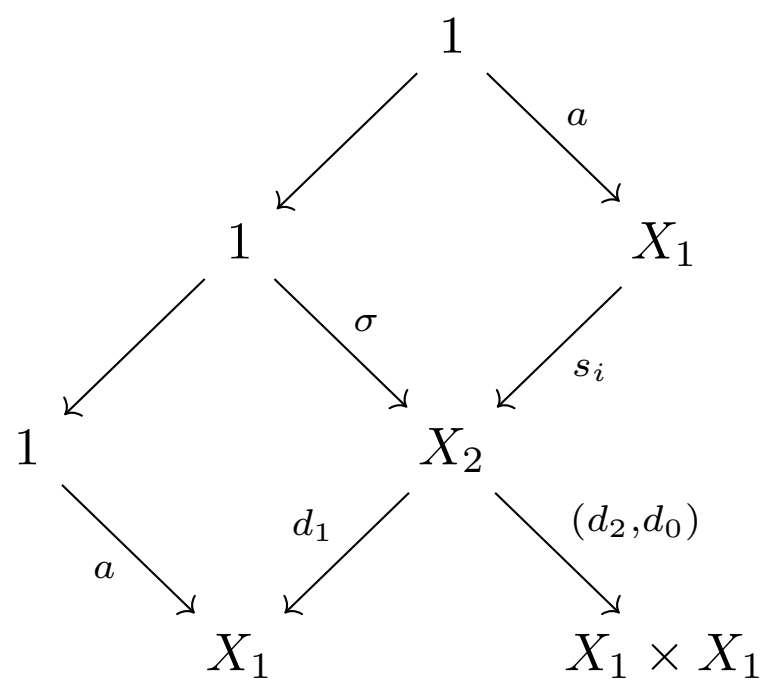
\[
  \begin{tikzcd}[column sep={4em,between origins}, row sep={9ex,between origins}]
   & & 1 \ar[ld]\ar[rd, "\xora"]& &  \\
  & 1 \ar[ld] \ar[rd, "\sigma"] && X_1 \ar[ld, "s_i"] 
  & \\
  1  \ar[rd, "\xora"'] && X_2 \ar[ld, "d_1"'] \ar[rd, "{(d_2,d_0)}"] && 
  \\
  & X_1 && X_1\times X_1 &
  \end{tikzcd}
    \]
  Here the composite rectangle is a pullback since $s_i$ is a section to $d_1$.
  Therefore the top square is a pullback too. But this shows that $\sigma =
  s_0(\xora) = s_1(\xora)$, and therefore $d_2(\sigma) = d_0(\sigma) = \xora$, so indeed we
  have $\Delta(\xora)=\xora\tensor \xora$.
\end{remark}

\begin{blanko}[Indecomposable identities]\label{indecomposable}
  Recall from \cite{Lawvere-Menni}, \cite[5.5]{GKT2} that a decomposition space has {\em indecomposable identities}  when
  we have a pullback
  \[
    \begin{tikzcd}
  X_0 \drpullback \ar[d,"="'] \ar[r, "s_0 s_0"] & X_2 \ar[d, "d_1"]  \\
  X_0 \ar[r, "s_0"'] & X_1  .
  \end{tikzcd}
\]
 In words, if the long edge of a triangle is degenerate, then the whole 
 triangle is doubly degenerate.
 When $X$ is the nerve of a category, it means that only identities can compose 
 to an identity.
\end{blanko}

\normalsize
The following is now clear from the definitions:
\begin{lemma}\label{degen are grouplike}
Let $X$ be a decomposition space. Then $s_0 \colon X_0 \to X_1$ is group-like
if and only if
  $X$ is complete and has indecomposable identities.
\end{lemma}

  Recall that a simplicial map is 
  {\em conservative} when it is cartesian on degeneracy maps. (For simplicial maps 
  between decomposition spaces, this can be measured on the first $s_0$.)
By the pasting law for pullbacks, we have:

\begin{lemma}\label{cons/complete+indecomp}
Suppose $Y\to X$ is conservative. 
Then $Y$ is complete if $X$ is, and $Y$ has indecomposable identities if $X$ does. \qed
\end{lemma}

Hitherto, in the theory of decomposition spaces, the gradings considered have
been only the length grading, namely the case where the canonical length
filtration is actually a grading. The filtration length of a $1$-simplex $\xora\in
X_1$ in a complete decomposition space $X$ is the maximal dimension of a 
nondegenerate simplex
with long edge $\xora$. The existence of the length filtration on $X$
is essentially the
Möbius condition~\cite{GKT2}, ensuring a Möbius
inversion principle. For the present purposes it is essential to go beyond length
gradings, because for our main object of study, 
the coalgebra of quasisymmetric functions, the relevant grading
is {\em not} the length 
grading, but rather the total grading, cf.~\ref{ex:Q-grading}.

\begin{blanko}[Grading of a decomposition space]\label{grading}
  We define a {\em grading} of a decomposition space $X$ to be a simplicial map to
  $B\N$. This means that every $1$-simplex $a\in X_1$ has associated to it a
  natural number and that for any higher simplex, say a $2$-simplex $\sigma$, we
  have
  \begin{equation*}
  \deg ( d_1(\sigma) ) = \deg( d_2(\sigma)) + \deg(d_0(\sigma)) .
  \end{equation*}
  as well as
  \begin{equation*}
  \deg( s_0(x) ) = 0
  \end{equation*}
  for any $0$-simplex $x \in X_0$.

  The grading defines a disjoint-union decomposition of $X_1$ by degree,
  $$
  X_1 = \sum_{k \in \N} X_1^{(k)} .
  $$
  We are particularly interested in the zeroth graded piece, which is thus the 
  pullback
  \begin{equation}\label{X^(0)}
  \begin{tikzcd}
  X_1^{(0)} \drpullback \ar[d] \ar[r] & X_1 \ar[d] \\
  1 \ar[r, "\name{0}"'] & \N  .
  \end{tikzcd}
  \end{equation}
  
  \end{blanko}

To say that a grading  $X\to B\N$ is conservative is to say that we have the
pullback square
\begin{equation}
  \label{eq:cons}
\begin{tikzcd}
  X_0 \drpullback \ar[d] \ar[r, "s_0"] & X_1 \ar[d] \\
  1 \ar[r, "\name{0}"'] & \N .
  \end{tikzcd}
\end{equation}
Comparing the two pullback squares, we find:

  \begin{lemma}\label{0=0}
    A grading is conservative precisely 
    when the canonical map $X_0 \to X_1^{(0)}$ is an equivalence.
  \end{lemma}

\begin{blanko}[Proto-connected graded decomposition spaces]\label{proto}
  We say that a graded decomposition space is {\em \proto-connected} when 
  $X_1^{(0)} \to X_1$ is group-like. (This is the decomposition-space analog of
  having degree zero spanned by group-like elements.)
\end{blanko}

\begin{proposition}\label{proto-connected = conservative}
  A graded decomposition space $X \to B\N$ is \proto-connected if and only if
  the grading is conservative.
\end{proposition}

\begin{proof}
  Suppose first that $X\to B\N$ is conservative.
  Since $B\N$ is complete and has indecomposable identities, 
  also $X$ is complete and has indecomposable identities (by  \cref{cons/complete+indecomp}), and then \cref{degen are grouplike}
  tells us 
  that $X_0 \to X_1$ is group-like. But to say that $X\to B\N$ is conservative is
  to say that $X_0 \simeq X_1^{(0)}$ (by \cref{0=0}). We therefore conclude that $X_1^{(0)}\to 
  X_1$ is group-like, which is the definition of being \proto-connected.
  
  Conversely, if $X$ is \proto-connected,
  then {the pullback square~\eqref{X^(0)} shows that} the map $X_1^{(0)} \to X_1$ is mono.
  It follows from the pullbacks in the diagram
  \[
  \begin{tikzcd}
  X_0 \ar[rrd, bend left, "\id"] \ar[rddd, bend right] \ar[rd, dotted] & &  \\
  & X_1^{(0)} \drpullback \ar[r] \ar[d, "="'] & X_0 \ar[d] \\
  & X_1^{(0)} \drpullback \ar[r] \ar[d] & X_1 \ar[d] \\
  & 1 \ar[r, "\name 0"'] & \N 
  \end{tikzcd}
  \]
  that also $X_1^{(0)} \to X_0$ is mono.
  This diagram shows that $X_1^{(0)} \to X_0$ is furthermore a retraction, 
  and since it is also a monomorphism, it is therefore an
  equivalence. But by \cref{0=0}, to have $X_0 \simeq X_1^{(0)}$ is precisely the condition that $X\to B\N$ be conservative.
\end{proof}

\begin{example}\label{ex:Q-grading}
  The {\em total grading} (the grading relevant to the universal property) on the decomposition space of quasisymmetric functions
  $Q = J(B\N)$  is precisely the
  canonical map $B\N \leftarrow J(B\N)$ from~\ref{J}, which is conservative by inspection or by \cref{counit conservative}.
\end{example}

We finish this section with a few remarks on the relationship between gradings and the
length filtration.
If $X$ admits the length filtration, then by evaluating $\ell_X \colon X_1 \to \N$ on the principal edges, we have a map $X_n \to X_1^{\times n} \to \N^{\times n} = B\N_n$ for each $n$.
We are interested in the situation when this constitutes a simplicial map $X\to B\N$, i.e.\ when the filtration is a grading.

\begin{lemma}
If a complete decomposition space $X$ admits the length filtration and this is a grading, then the grading is conservative.
\end{lemma}
\begin{proof}
  The filtration assigns a natural number to each $1$-simplex, so we have a map
  $X_1\to \N$ and can ask for the pullback condition
  \[
  \begin{tikzcd}
    X_0 \drpullback \ar[d] \ar[r, "s_0"] & X_1 \ar[d] \\
    1 \ar[r, "\name{0}"'] & \N
    \end{tikzcd}
  \] 
  even without knowing that it is part of a simplicial map. The length filtration
  implies that $X$ has indecomposable
  identities~\cite{GKT2}, and now it is easy to see
  that the pullback condition is verified: 
  a degenerate $1$-simplex cannot be the long edge of any higher-dimensional nondegenerate  
  simplex because of the indecomposable-identities condition, and a non-degenerate
  $1$-simplex is clearly of length at least $1$. This establishes the pullback
  without even knowing whether the filtration is a grading. If it 
  is a grading, then the condition says it is conservative.
\end{proof}

We also have implications of a grading for the length filtration, derived from the 
following lemma (from \cite[Lemma 5.4.2]{GKT:Crapo}).

\begin{lemma}\label{lem:locally finite length}
  Let $f:Y\to X$ be a conservative map between decomposition spaces.
  For every $\yora\in Y_1$ we have the inequality of lengths
  $$
  \ell_Y(\yora) \leq \ell_X(f(\yora)).
  $$
  In particular, if $X$ is (complete and) of locally finite length, then so is $Y$.
\end{lemma}

Since $B\N$ is of locally finite length, we get:

\begin{corollary}
  A (conservatively) graded 
  decomposition space is automatically of locally finite length
  (but the grading does not have to coincide with the length filtration).
\end{corollary}

\begin{example}
  The (conservative) total grading $Q \to B\N$ of \cref{ex:Q-grading}, sends, for example, $(216)$ to $9$. Both these decomposition spaces are of locally
  finite length. The inequality of \cref{lem:locally finite length} is illustrated by
  $\ell_Q(216) = 3$ while $\ell_{B\N}(9)=9$.
\end{example}
  
\begin{remark}
  One may wonder if it is interesting to ask a grading to be culf rather than just 
  conservative. It is interesting, but very restrictive: it is shown in 
  \cite{Hackney-Kock:2210.11192} that a decomposition space admits a culf 
  functor to $B\N$ if and only if it is a free decomposition space.
\end{remark}

\begin{blanko}[More general notions of grading]
  In this section we have studied gradings with values in $B\N$ and shown that
  conservativity of this map implies the desired property (\proto-connectedness) that
  degree zero is spanned by group-like elements. For the proof we used only that $B\N$
  is complete and has indecomposable identities. For the main theorem in
  Section~\ref{sec main thm proof}, we work with more general recipients of the
  grading, namely a complete decomposition space $M$ with indecomposable identities,
  and a conservative map $X\to M$. 
\end{blanko}

\section{Zeta functions}
\label{sec:zeta}

The main point of this section is to make explicit the notions of standard and 
special zeta function, and to explain the notion of grounded functional (included 
in the definition of zeta function).

\begin{blanko}[Zeta functions in algebraic combinatorics]
  Linear functionals on a combinatorial coalgebra form the convolution
  algebra. If the coalgebra is an incidence coalgebra (of a decomposition space)
  the convolution algebra is the {\em incidence algebra}. For many purposes it is
  necessary to impose the condition that a functional $\phi$ should take value $1$ on
  group-like elements. This means that $\phi$ agrees with $\epsilon$ on group-like
  elements.
  This condition ensures that the functional has a convolution
  inverse, constructed as a geometric series
  \begin{equation}\label{conv-inv}
  \phi^{*-1} = \sum_{k\in\N} (-1)^k (\phi-\epsilon)^{*k} .
  \end{equation}

  In the present context, we are not concerned 
  with any inversion, but the condition is simply necessary
  in order for the universal property to hold.%
        \footnote{The condition that the functional takes value 1 on group-like elements was missed in the published version of \cite{Aguiar-Bergeron-Sottile}, 
        because it is automatic in the
        bialgebra case (where the zeta function is required to be a character).
        The necessity of imposing this extra assumption was first pointed out by
        Darij Grinberg. A corrected version of the paper is available from
        Aguiar's homepage \url{https://pi.math.cornell.edu/~maguiar/pubs.html}.}
  (A counter example can be given by the zero functional, which cannot be
  induced by the special zeta function of $\QSym$.)
  We follow Aguiar--Bergeron--Sottile~\cite{Aguiar-Bergeron-Sottile} in calling
  linear functionals satisfying this condition {\em zeta functions}. In the case
  of incidence algebras, the word zeta function is often used in a more
  restricted sense, which is the following.
\end{blanko}

\begin{blanko}[Standard zeta functions]
  The {\em standard zeta function} on a combinatorial coalgebra is the zeta
  function that takes all basis elements to $1$. The Möbius function is by
  definition the convolution inverse of the standard zeta
  function~\cite{Rota:Moebius}.
\end{blanko}

\begin{blanko}[The zeta function of $\QSym$ and special zeta functions]
  The zeta function on $\QSym$ featured in the universal property is not the 
  standard zeta function but rather the
  {\em special} zeta function (cf.\ \ref{Q-gradings-and-special-zeta}). 
  What the word ``special'' should mean in
  general, is not obvious at the level of classical algebraic
  combinatorics, but we shall see at the objective level that the decomposition
  space whose incidence coalgebra is $\QSym$ is {\em free}, and that the special zeta
  function is free on a standard zeta function. So in a moment we shall define
  special zeta functions as those that arise in that way.
\end{blanko}

\begin{blanko}[Zeta functions at the level of decomposition spaces]\label{zeta decomp}
  At the objective level where decomposition spaces $X$ play the role of coalgebras,
  a linear functional is a span
  $$
  X_1 \leftarrow F \to 1 .
  $$
  In order to give an objective formulation of the zeta-function condition and
  make sense of the subtraction $\phi-\epsilon$ in the inversion formula 
  \eqref{conv-inv}, the space $F$ 
  should receive a map from $X_0$ factoring $s_0$, and that this map should be 
  a monomorphism, as in the diagram
\begin{equation}\label{tilde-diagram}
\begin{tikzcd}
 & F \ar[d] && \\
X_0 \ar[ru, rightarrowtail, shorten >=1.2mm, shorten <=2.0mm] \ar[r, shorten >=1.4mm, shorten <=2.0mm] 
& X_1 
\ar[l, shift left=1.5]
\ar[l, shift right=1.5]
\ar[r, shift left=1.3, shorten >=1.4mm, shorten <=2.0mm]
\ar[r, shift right=1.3, shorten >=1.4mm, shorten <=2.0mm]
& X_2 
\ar[l, shift left=2.6]
\ar[l] 
\ar[l, shift right=2.6] \rar[dotted, no head]
& {} .
\end{tikzcd}
\end{equation}
If $X$ is complete, the zeta function condition says the following is a pullback
  \[
  \begin{tikzcd}
  X_0 \ar[d, equal] \drpullback \ar[r]  & F \ar[d]  \\
  X_0 \ar[r, "s_0"'] & X_1 .
  \end{tikzcd}
  \]
We call such $F$ {\em grounded}. We define a {\em zeta function} to be a grounded linear 
functional.
The picture \eqref{tilde-diagram} above leads to the category $\Tsimpl$ 
introduced in the next section, which is the setting for most of the technical
results of the paper.
\end{blanko}

\begin{blanko}[Standard zeta function]
  At the level of decomposition spaces, the {\em standard} zeta function is simply the
  span $X_1 \stackrel{=}\leftarrow X_1 \to 1$, which with the canonical 
  pointing $s_0: X_0\to X_1$ is clearly grounded. 
  Rather than saying anything about an identity map, we shall say that a functional
  $X_1 \leftarrow F$ is {\em standard} if this map is an equivalence.
  See also \cref{section standard}.
\end{blanko}

\begin{blanko}[Functorialities]
  \Cref{intro theorem A} involves the statement that a given zeta
  function $F$ (on some complete decomposition space $X$) is induced from the special
  zeta function on the decomposition space $Q$ along a certain ikeo-culf span. Let us
  briefly explain what it means to induce zeta functions. Coalgebra homomorphisms act
  on zeta functions by precomposition. At the objective level this looks different
  for culf and for ikeo maps. If $Y \to X$ is a culf map, then a zeta function
  $X_1\leftarrow F \to 1$ induces the zeta function given by pullback, $Y_1
  \leftarrow Y_1\times_{X_1} F \to 1$. 
  As an example, the standard zeta function 
  $X_1
  \stackrel{=}\leftarrow X_1 \to 1$ induces the standard zeta function $Y_1
  \stackrel{=}\leftarrow Y_1 \to 1$. For ikeo maps it is instead about
  postcomposition: if $Y \to X$ is ikeo, then a functional $Y_1\leftarrow F \to
  1$ induces the functional given as the composite $X_1 \leftarrow
  Y_1\leftarrow F \to 1$. 
  Groundedness is preserved by \emph{conservative} ikeo maps.
\end{blanko}

\begin{blanko}[Special zeta function of $Q=J(B\N)$]
  The special zeta function of $Q$ is given by the span
  $$
  Q_1 \longleftarrow \{\varnothing\} + \N_+ \longrightarrow 1 .
  $$
  Recall that $Q_1= (\N_+)^*$ is the set of all words in the alphabet $\N_+$;
  the set in the middle of the span consists of the empty word and all the
  one-letter words, expressing the fact that the zeta function of $\QSym$
  considered by Aguiar--Bergeron--Sottile~\cite{Aguiar-Bergeron-Sottile} sends
  the quasisymmetric monomials $M_{\{\}}$ and $M_n$ to $1$ and all other basis
  elements to $0$. The universal property of $\QSym$ refers to this choice of
  zeta function, stating that $(\QSym, \zeta_{\operatorname{special}})$ is the
  terminal object in the category of (connected) graded coalgebras equipped with
  a zeta function.
\end{blanko}

\begin{blanko}[Example (standard zeta function induced from special)]
  It was  noted in \cite{GKT:QSym} that
  the standard zeta function on a complete decomposition space $X$ is induced from the
  special zeta function on $Q$ along the span
  $$
  X \longleftarrow J(X) \longrightarrow J(B\N) = Q
  $$
  (but no attempt was made to establish the universal property, and the case of 
  standard zeta function is also very particular). With the observation from
  \cite{Hackney-Kock:2210.11192} that $J$ is an instance of the free 
  construction,
  $J(X) = \lilj_! (\nondeg X)$ (cf.~\ref{J}), and given the notion of grading $\nu: X \to
  B\N$ (cf.~\ref{grading}), the right-pointing map is just $\lilj_!(\nu)$.
\end{blanko}

  The notion of special zeta function can be defined 
  uniformly for any free decomposition space: 

\begin{blanko}[Special zeta function of a free decomposition space]\label{zeta of free}
  For $Y\coloneq \lilj_!(A)$ the free decomposition space on $A\in \pre(\inert)$, we
  define the {\em special zeta function} to be the span
  $$
  Y_1 \leftarrow A_0+A_1 \to 1 .
  $$
  Recall that $Y_1 = \sum_n A_n$, so the leftward map  in the span is just the 
  inclusion of the first two summands. 
  
  In the case $A= \nondeg X$ (cf.\ \ref{J}), 
  the special zeta function on $J(X) = \lilj_!(\nondeg X)$ is thus
  $$
  \sum_n \nondeg X_n \longleftarrow \nondeg X_0+ \nondeg X_1 \longrightarrow 1 .
  $$
  Since the composite 
  $$
  X_1 \leftarrow \sum_n \nondeg X_n \leftarrow \nondeg X_0+ \nondeg X_1
  $$
  is the identity (expressing the sum splitting of $X_1$ into degenerate and
  nondegenerate simplices), we see that the coalgebra homomorphism given by the
  ikeo map $X \leftarrow J(X)$ sends the special zeta function
  of $J(X)$ to the standard zeta function of $X$. 
  In slogan form: {\em ``special means free on standard''}.
\end{blanko}

\section{The Expanded Simplex Category}\label{sec shape categories}

In order to encode zeta functions for decomposition spaces we expand the simplicial category $\simpl$ to $\Tsimpl$, incorporating the zeta function at the combinatorial level, working with presheaves on $\Tsimpl$ rather than ordinary simplicial objects.
The trouble of setting up this category $\Tsimpl$ pays off: many of the conditions we need to handle in connection with zeta functions are absorbed into well-appreciated notions such as conservative, culf and ikeo.
Conceptually this reflects the philosophy of Aguiar--Bergeron--Sottile~\cite{Aguiar-Bergeron-Sottile} that ``combinatorial coalgebra'' means ``(graded) coalgebra with a zeta function'' -- we treat it as part of the structure of our objects, rather than something added on.

To expand the simplex category $\simpl$, we factor the
codegeneracy map $[0] \leftarrow [1]$ through a new object $\Tob$.
We write $\Tsimpl \coloneq \simpl \amalg_{[1]} [2]$ for this new category,\footnote{This pushout may be taken either in $\cat$ or $\catinf$ -- see \cref{app subdividing}.} and $\sigr \colon [1] \to
\Tob$ and $\sigl \colon \Tob \to [0]$ for its new generators\footnote{This category $\Tsimpl$ previously appeared in \cite[Remark 1.3]{GagnaHarpazLanari:OEAMI2C} under the name $\simpl^+$ (this is a full subcategory of the category $t\simpl$ from \cite[Section 1.1]{OzornovaRovelli:MSINCPSSSPSS} or \cite[Observation 12]{Verity:WCS1BHT}).}
(with $\sigl \sigr = \sigma^0$).
So $\simpl \subset \Tsimpl$ is a full subcategory
inclusion. The picture to keep in mind is:
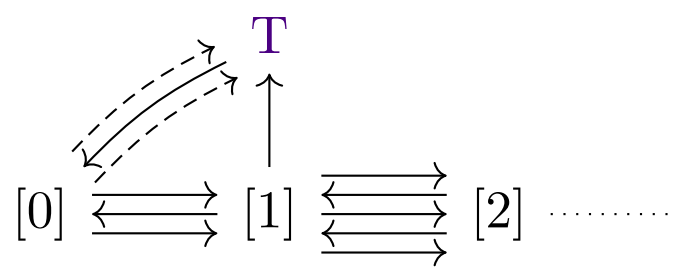
\[ \begin{tikzcd}
& \Tob \dlar[bend right=10] \\
{[0]} 
\rar[shift left=1.5] \rar[shift right=1.5] 
\urar[shift left=1.5, bend left=10, dashed] \urar[shift right=1.5, bend left=10, dashed] 
& {[1]} \lar \uar 
\rar
\rar[shift right=3]
\rar[shift left=3]
& {[2]} \rar[dotted,no head] 
\lar[shift right=1.5]
\lar[shift left=1.5]
& {}
\end{tikzcd} \]
The following list contains all morphisms of $\Tsimpl$ (this is easily checked directly; alternatively, it follows from \cref{lem pushout calc}):
\begin{itemize}
\item morphisms in $\simpl$
\item $\hom(\Tob,[n]) \cong \hom([0],[n])$ via $\alpha \circ \sigl \mapsfrom \alpha$ (so for a fixed $n$ there are $n+1$ maps $\Tob\to [n]$)
\item $\hom([n],\Tob) \cong \hom([n],[1])$ via $\sigr \circ \alpha \mapsfrom \alpha$
\item $\hom(\Tob,\Tob) \cong \hom([0],[1]) \amalg \{\id_\Tob \}$ via $\sigr \circ \delta^i \circ \sigl \mapsfrom \delta^i$; this is a 3-element set.
\end{itemize}

There is a retraction $\lilp \colon \Tsimpl \to \simpl$ sending $\Tob$ to $[1]$, $\sigr$ to $\id_{[1]}$, and $\sigl$ to $\sigma^0$.
The retraction $\lilp$ is right adjoint to the inclusion $\lili \colon \simpl \hookrightarrow \Tsimpl$, with counit $\lili \lilp \Rightarrow \id_{\Tsimpl}$ given at $\Tob$ by $\sigr$, so $\simpl$ is a coreflective subcategory of $\Tsimpl$.

\begin{blanko}[The active-inert factorization system]
The simplicial category $\simpl$ has an active-inert factorization system $(\simpl^\actrm, \simpl^\intrm)$ where active maps are endpoint-preserving, and inert maps are distance-preserving.
This factorization system extends to a factorization system on $\Tsimpl$, 
by specifying that
\begin{itemize}
\item

both $\sigr$ and $\sigl$ are active, 

\item the maps $\sigr \circ \delta^j \colon [0] \to \Tob$ for $j=0,1$ are inert.

\end{itemize}
For general maps in $\Tsimpl \setminus \simpl$:
\begin{itemize}
\item $\id_\Tob$: both active and inert.
\item $\alpha \circ \sigl$: since $\alpha$ is a map with domain $[0]$, it is automatically inert. This composite is active if and only if $\alpha$ is an identity, otherwise it is neither active nor inert.
\item $\sigr \circ \alpha$: this map is active if $\alpha$ is active. It is inert if and only if $\alpha$ is a non-identity inert map (namely $\delta^0,\delta^1 \colon [0] \rint [1]$).
\item $\sigr \circ \delta^i \circ \sigl$: These maps are neither active nor inert. 
The active-inert factorization is $(\sigr \delta^i) \sigl$.
\end{itemize}
\end{blanko}

\begin{blanko}[The negative-positive factorization system]
Notice that $\sigr$ is a monomorphism and $\sigl$ is a retraction/split-epimorphism (since
$\id_{[0]} = \sigma^0 \delta^0 = \sigl (\sigr \delta^0)$). One can check that
$\sigl$ is the \emph{only} non-identity retraction of
$\Tsimpl$ which is not contained in $\simpl$. Each monomorphism in $\simpl$
remains a monomorphism in $\Tsimpl$. Maps of the form $\alpha \circ \sigl$ or
$\sigr \circ \delta^i \circ \sigl$ are not monomorphisms, as witnessed by the
distinct endomorphisms $\sigr \circ \delta^0 \circ \sigl$ and $\sigr \circ
\delta^1 \circ \sigl$ of $\Tob$. Finally, a map of the form $\sigr \circ \alpha
\colon [n] \to \Tob$ is a monomorphism if and only if $\alpha \colon [n] \to
[1]$ is a monomorphism, i.e.\ just when $\alpha$ is the identity or $\delta^i
\colon [0] \to [1]$. We observe that there is a retraction--mono factorization
system $(\Tsimpl^-,\Tsimpl^+)$; in particular we have $\Tsimpl^- \cap \simpl =
\simpl^-$ is the subcategory of surjections and $\Tsimpl^+ \cap \simpl =
\simpl^+$ is the subcategory of injections. Notice that $\Tsimpl^- \subset
\Tsimpl^\actrm$ and $\Tsimpl^\intrm \subset \Tsimpl^+$, so in fact we have a
ternary factorization system.
\Cref{figure generating morphisms} below contains a summary for the generating morphisms of $\Tsimpl$, as well as other categories defined in what follows.
\end{blanko}

\begin{blanko}[Contractive maps]
An order-preserving map $f\colon [k] \to [n]$ is \emph{contractive} if it is
weakly distance-decreasing, i.e.\ if $f$ satisfies $f(i+1) \leq f(i)+1$ for all
$i=0,\dots,k-1$. We let $\contr \subset \simpl$ be the wide subcategory of
contractive maps. This is the subcategory of $\simpl$ spanned by surjections and
inert maps.

The category $\Tcontr \subset \Tsimpl$ is obtained from $\contr$ by factoring $\sigma^0 \colon [1] \to [0]$ in $\contr$ through a new object $\Tob$. It is the
subcategory of $\Tsimpl$ spanned by the inert maps $\Tsimpl^\intrm$, the
retractions $\Tsimpl^-$, and $\sigr \colon [1] \to \Tob$ (notice $\sigl \in
\Tsimpl^- \subset \Tcontr$).
The category $\Tcontr$ admits an active-inert factorization system induced from 
the active-inert factorization system on $\Tsimpl$: the active maps are $\Tcontr^\actrm = \Tsimpl^\actrm \cap \Tcontr$ (generated by $\Tsimpl^-$ and $\sigr$) and the inert maps are 
$\Tcontr^\intrm = \Tsimpl^\intrm \cap \Tcontr = \Tsimpl^\intrm$.
\end{blanko}

The \emph{monomorphisms} in $\Tcontr$ form another important category which
we call $\Tinert$ in this paper. Let us quickly characterize the maps that
appear in this subcategory: it is immediate that $\Tcontr \cap \Tsimpl^+
\subseteq \Tinert$, and we will show that the reverse inclusion holds.
First, if $[p] \to [q]$ in $\Tcontr$ is mono, then it is injective as a map
in $\simpl$ (since the maps $[0] \rint [p]$ are inert and hence in
$\contr$), hence is in $\simpl^+$. Mimicking the reasoning in $\Tsimpl$,
the only other monomorphisms in $\Tcontr$ are $\id_\Tob$, $\sigr$, and
$\sigr \circ \delta^i \colon [0] \to \Tob$. So $\Tinert$ takes the
following form.
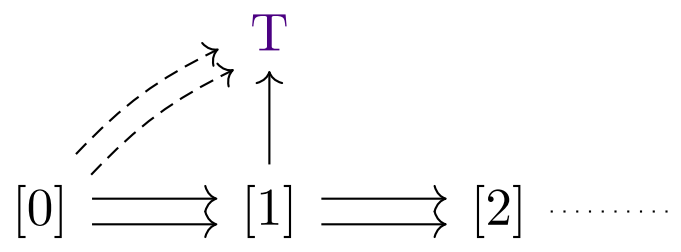
\[ \begin{tikzcd}
& \Tob  
\\
{[0]} 
\rar[shift left=1] \rar[shift right=1] 
\urar[shift left=1, bend left=10, dashed] \urar[shift right=1, bend left=10, dashed] 
& {[1]}
\uar 
\rar[shift right=1]
\rar[shift left=1]
& {[2]} \rar[dotted,no head] 
& {}
\end{tikzcd} \]
We notice that $\Tinert \cap \simpl = \simpl^\intrm$ is the subcategory of
inert maps, and we also write $\inert = \simpl^\intrm$ for consistency
between the three situations. Notice that the active-inert factorization system on $\Tsimpl$ restricts to a factorization system $(\Tinert^\actrm, \Tinert^\intrm)$ on $\Tinert$ where
$\Tinert^\actrm$ has a single nonidentity map $\sigr$, and $\Tinert^\intrm
= \inert \sqcup \{ \id_\Tob, \sigr \delta^0, \sigr \delta^1 \}$.

\begin{notation}\label{notation for functors}
We write $\lilk \colon \Tcontr \to \Tsimpl$ and $\lilh \colon \Tinert \to \Tcontr$ for the inclusions, and $\lilj \coloneq \lilk \circ \lilh$.
\end{notation}
We give restrictions of functors the same names, as on the bottom line below:
\begin{equation}
\label{commutative ipjlk}
\begin{tikzcd}[column sep={45pt,between origins}, row sep={45pt,between origins}]
\Tinert \rar{\lilh} 
\ar[d, shift left=1.6, pos=0.54, "\lilp"] \ar[d, phantom, "\dashv"{font=\footnotesize}]
& \Tcontr \rar{\lilk} 
\ar[d, shift left=1.6, pos=0.54, "\lilp"] \ar[d, phantom, "\dashv"{font=\footnotesize}]
& \Tsimpl
\ar[d, shift left=1.6, pos=0.54, "\lilp"] \ar[d, phantom, "\dashv"{font=\footnotesize}]
\\
\inert \ar[u, shift left=1.6, "\lili"] \ar[r, "\lilh"'] 
& \contr \ar[u, shift left=1.6, "\lili"] \ar[r, "\lilk"'] 
& \simpl \ar[u, shift left=1.6, "\lili"]
\end{tikzcd}  
\end{equation} 

\begin{blanko}[Restricted ternary factorization systems]\label{ss ternary}
The category $\Tsimpl$ has a ternary factorization system, and the retraction--mono factorization system $(\Tsimpl^-,\Tsimpl^+)$ also induces factorization systems on $\Tcontr$ and $\Tinert$.
We have the following negative-positive factorization systems which play a role below: 
\[ 
\begin{gathered}
(\Tcontr^-,\Tcontr^+) = (\Tsimpl^- \cap \Tcontr, \Tsimpl^+ \cap \Tcontr) = (\Tsimpl^-, \Tinert)  \\
(\Tinert^-,\Tinert^+) = (\Tsimpl^- \cap \Tinert, \Tsimpl^+ \cap \Tinert) = (\ob(\Tsimpl), \Tinert).
\end{gathered}
\]
\end{blanko}

\begin{figure}[htb]
\begin{tikzpicture}
\begin{scope}[shift={(0,-3.2)}]
\node at (2.1, 2.4) {$\overbrace{\phantom{xxxxxxxxxxxxxxxxxx}}$};
\node at (2.1, 2.9) {\text{active}};
\node at (5.2, 2.4) {$\overbrace{\phantom{xxxxxxx}}$};
\node at (5.2, 2.9) {\text{inert}};

\node at (4.2, -0.3) {$\underbrace{\phantom{xxxxxxxxxxxxxxxxxxx}}$};
\node at (4.2,-0.8) {\text{positive}};
\node at (1.0, -0.3) {$\underbrace{\phantom{xxxxxxx}}$};
\node at (1.0,-0.8) {\text{negative}};

\draw (0,0) -- (6.5,0) -- (6.5,2.1) -- (0,2.1) -- (0,0);
\node at (-1,1.6) {$\Tsimpl$ :};
\node at (1.0,1.6) {$\sigma^i, \sigl$};
\node at (3.0,1.6) {$\sigr,\  \delta^j$};
\node at (5.2,1.6) {$\delta^\bot, \delta^\top$};
\node at (-1,1) {$\Tcontr$ :};
\node at (1.0,1) {$\sigma^i, \sigl$};
\node at (2.7,1) {$\sigr$};
\node at (5.2,1) {$\delta^\bot, \delta^\top$};
\node at (-1,0.4) {$\Tinert$ :};
\node at (2.7,0.4) {$\sigr$};
\node at (5.2,0.4) {$\delta^\bot, \delta^\top$};
\end{scope}
\end{tikzpicture}
\caption{Locations of generating morphisms of $\Tsimpl$ ($i$ arbitrary, $j\neq \top, \bot$)}\label{figure generating morphisms}
\end{figure}

\section{The 2-Segal condition and the sturdy condition}
\label{sec object prop}

In this section we give two important exactness conditions for presheaves which will be used throughout. 

\begin{convention}
In what follows we work with presheaves in a fixed $\infty$-category $\basecat$ satisfying the conditions in \ref{base category}, and $\pre(\mathcal{C})$ will be shorthand for $\fun(\mathcal{C}^\op, \basecat)$. 
This will continue until \cref{sec finiteness}, at which point we will specialize again to the case $\basecat = \spaces$.
\end{convention}

\begin{definition}
We say a presheaf $X\in \pre(\Tsimpl)$ is \emph{2-Segal} if each active-inert pushout square that exists in $\Tsimpl$ is sent to a pullback square.
\end{definition}
This amounts to asking that $\lili^* X \in \pre(\simpl)$ is 2-Segal \cite[\S2.3]{DyckerhoffKapranov:HSS} (or satisfies the decomposition-space axioms \cite[\S3]{GKT1}), which is a consequence of the following lemma.

\begin{lemma}\label{lem no more ai pushouts}
Each nontrivial active-inert pushout square in $\Tsimpl$ is in $\simpl$.
\end{lemma}
\begin{proof}
As $\lili \colon \simpl \to \Tsimpl$ is a left adjoint, each pushout in $\simpl$ is also a pushout in $\Tsimpl$.
The nonidentity inert maps in $\Tsimpl \setminus \simpl$ are $\sigr \delta^0, \sigr \delta^1 \colon [0] \to \Tob$, but no nonidentity active map has domain $[0]$, so there are no pushouts to consider with one leg in $\Tsimpl^\intrm \setminus \simpl$.
The nonidentity active maps in $\Tsimpl\setminus \simpl$ are $\sigl$ and $\sigr \circ \alpha$ with $\alpha \colon [n] \to [1]$ active in $\simpl$.
We may disregard $\sigl$, since there are no nonidentity inert maps with domain
$\Tob$.
Since the pushout of any inert map in $\simpl$ along $\alpha \in \simpl^\actrm$ exists and is again inert, the only case we must consider is $\sigr \in \Tsimpl^\actrm$. 

Any inert map with domain $[1]$ is in $\simpl$, so it suffices to show that for each nonidentity inert map $\iota \colon [1] \rint [m]$, the span
\[ \begin{tikzcd}
{[1]} \dar[-act,"\sigr"'] \rar[tail,"\iota"] & {[m]} \\
\Tob
\end{tikzcd} \]
does not admit a pushout.
Suppose it did: it would take one of the forms below
\[ \begin{tikzcd}
{[1]} \dar[-act,"\sigr"'] \rar[tail,"\iota"] & {[m]} \dar{\sigr \alpha}
&
{[1]} \dar[-act,"\sigr"'] \rar[tail,"\iota"] & {[m]} \dar{\gamma}
\\
\Tob \rar{?} & \Tob 
&
\Tob \rar{\beta \sigl} & {[p]}
\end{tikzcd} \]
for some $\alpha \colon [m] \to [1]$ or $\beta \colon [0] \to [p]$.

We first exclude the possibility of the left square.
Letting $i = \iota(0)$, the composite $\sigma^i \iota \colon [1] \to [m-1]$ is constant at $i$, and thus factors uniquely as $\name{i} \sigma^0 = \name{i} \sigl \sigr$ through $\Tob$.
There does not exist a $\chi \sigl \colon \Tob \to [m-1]$ satisfying $(\chi 
\sigl) (\sigr \alpha) = \sigma^i$, as the left-hand side is constant and the 
right-hand side is not (since $m > 1$).

To exclude the possibility of the right square, let $\rho \colon [m] \to [1]$ be the unique map with $\rho \iota = \id_{[1]}$, and consider the maps $\sigr \rho \colon [m] \to \Tob$ and $\id_\Tob \colon \Tob \to \Tob$.
We have $(\sigr \rho) \circ \iota = \sigr = \id_\Tob \circ \sigr$, but 
there is no map $\sigr \chi \colon [p] \to [1] \to \Tob$ that composes with $\beta \sigl$ to give the identity.
\end{proof}

We now want to approach the \emph{sturdy} presheaves (\cref{def new sturdy}), but we first need some preliminaries.
Recall that in $\simpl$, surjective maps are stable under pushout.
We wish to establish that negative maps in $\Tsimpl$ and $\Tcontr$ are also stable under pushout.
The following cancellation lemma is key for establishing when certain pushouts in $\Tsimpl$ are also pushouts in $\Tcontr$.

\begin{lemma}\label{lem cancellation}
  If $g\circ\sigma \in \Tcontr$, where $g$ is a map in $\Tsimpl$ and $\sigma$ is a map
  in $\Tsimpl^-$, then $g\in \Tcontr$.
\end{lemma}
\begin{proof}
If $g$ and $\sigma$ are both in $\simpl$, then we wish to show $g(i+1) \leq g(i) + 1$.
Let $j$ be the maximal element of $\sigma^{-1}(i)$, which implies $j+1$ is the minimal element of $\sigma^{-1}(i+1)$.
Then
\[
  g(i+1) = g\sigma(j+1) \leq g\sigma(j) +1 = g(i) + 1,
\]
so $g$ is contractive.
If $\sigma$ is not in $\simpl$, then it is either $\sigl \colon \Tob \to [0]$ or $\id_\Tob$; in the latter case we have $g = g\id_\Tob = g\sigma \in \Tcontr$.
In the former case, every map in $\Tsimpl$ with domain $[0]$ is in $\Tcontr$, hence $g\in \Tcontr$.

We are left with the situation where $\sigma$ is in $\simpl$, but $g$ is not.
Then $g = \sigr \alpha$ for a unique $\alpha \colon [m] \to [1]$.
Since $g \sigma = \sigr \alpha \sigma$ is in $\Tcontr$, the map $\alpha \sigma$ is in $\contr$, hence $\alpha \in \contr$ by the first paragraph.
It follows that $g\in \Tcontr$.
\end{proof}

\begin{lemma}\label{lem pushout stability negative}
Suppose $\sigma \colon a \twoheadrightarrow c$ is in $\Tsimpl^- = \Tcontr^-$ and $f\colon a \to b$ is in $\Tsimpl$.
Then there is a pushout square
\begin{equation}\label{diag neg pushout}
\begin{tikzcd}
a \rar[two heads, "\sigma"] \dar[swap]{f} \drar[phantom, "\ulcorner"very near end] & c \dar{g} \\
b \rar[two heads, "\bar \sigma"'] & d
\end{tikzcd}
\end{equation}
whose bottom arrow is negative.
If $f\in \Tcontr$, then the pushout \eqref{diag neg pushout} in $\Tsimpl$ is also a pushout in $\Tcontr$.
\end{lemma}
\begin{proof}
We'll begin with the second statement.
Assume we have the pushout \eqref{diag neg pushout} in $\Tsimpl$ with $f\in \Tcontr$.
Then $g\in \Tcontr$ by \cref{lem cancellation}.
Now suppose $u\colon b \to z$ and $v\colon c \to z$ are maps in $\Tcontr$ with $uf = v\sigma$, and let $m\colon d \to z$ be the unique map in $\Tsimpl$ with $m\bar\sigma = u$ and $mg = v$.
Then $m$ is in $\Tcontr$ by \cref{lem cancellation}, so we conclude that \eqref{diag neg pushout} is a pushout in $\Tcontr$.

We turn to the first statement. 
If all of $a,b,c$ are in $\simpl$, then the corresponding pushout in $\simpl$ is also a pushout in $\Tsimpl$ since $\lili \colon \simpl \to \Tsimpl$ is a left adjoint.

In the general setting, we may restrict to the case when $f$ is a generating morphism of $\Tsimpl$.
There is a single nonidentity map in $\Tsimpl^- \setminus \simpl^-$, namely $\sigl \colon \Tob \to [0]$.
This is also the only generator of $\Tsimpl$ with domain $\Tob$, and the square
\[ \begin{tikzcd}
\Tob \rar{\sigl} \dar[swap]{\sigl} & {[0]} \dar{=}  \\
{[0]} \rar[swap]{=} & {[0]}
\end{tikzcd} \]
is a pushout since $\sigl$ is an epimorphism.
If $\sigma$ is in $\simpl$ but the generating morphism $f$ is not, then the only nontrivial situation possible is $f=\sigr \colon [1] \to \Tob$ and $\sigma \colon [1] \twoheadrightarrow [0]$.
We now observe that the square
\begin{equation}\label{diag grounded pushout}
\begin{tikzcd}
{[1]} \dar[swap]{\sigr} \rar[two heads,"\sigma^0"] & {[0]} \dar{=} \\
\Tob \rar{\sigl} & {[0]}
\end{tikzcd} 
\end{equation}
is a pushout.
This happens since any nonidentity map $u \colon \Tob \to z$ factors uniquely as $u = v\sigl$, while $\id_\Tob \circ \sigr = \sigr$ does not factor through $[0]$.
\end{proof}

\begin{definition}\label{def new sturdy}
A presheaf $A\in \pre(\Tcontr)$ is \emph{sturdy} if it sends each pushout square with one leg negative to a pullback.
A presheaf $X\in \pre(\Tsimpl)$ is sturdy just when $\lilk^*X \in \pre(\Tcontr)$ is sturdy.
\end{definition}

By looking at special cases for $f$ in \cref{lem pushout stability negative}, we arrive at the following alternative characterization of sturdy.

\begin{proposition}\label{prop sturdy char}
A presheaf $A\in \pre(\Tcontr)$ is sturdy if and only if the following conditions hold:
\begin{enumerate}
\item $A$ is \emph{grounded}, i.e.\ the square
\[ \begin{tikzcd}
A_0 \ar[d, "="'] \rar{\sigl^*} & A_\Tob \dar{\sigr^*} \\
A_0 \rar[swap]{s_0} & A_1
\end{tikzcd} \]
is a pullback. \label{sturdy grounded}
\item Each active-inert pushout that exists in $\Tcontr$ is sent to a pullback.\label{sturdy stiff}
\item Each negative map in $\Tcontr$ is sent to a monomorphism, i.e.~if $c \leftarrow a$ is negative, then $A_c \to A_c \times_{A_a} A_c$ is an equivalence.\label{sturdy complete}
\end{enumerate}
\end{proposition}
\begin{proof}
Certainly any sturdy presheaf satisfies the three conditions. 
The square in \eqref{sturdy grounded} is the image of the pushout \eqref{diag grounded pushout} in the proof of \cref{lem pushout stability negative}.
By the proof of \cref{lem no more ai pushouts}, active-inert pushouts in $\Tcontr$ are the same thing as negative-inert pushouts, so \eqref{sturdy stiff} holds. 
Finally, given a negative map $\sigma$, since $\sigma$ is epi the square below left is a pushout, hence sent to the pullback below right.
\[
\begin{tikzcd}
c \drpullback & c \lar[swap]{=} \\
c \uar["="] & a
\lar[two heads, "\sigma"] \uar[swap, two heads]{\sigma}
\end{tikzcd} 
\ \rightsquigarrow \
\begin{tikzcd}
A_c \rar{=} \dar[swap]{=} \drar[phantom, "\lrcorner"very near start] & A_c \dar{\sigma^*} \\
A_c \rar[swap]{\sigma^*} & A_a
\end{tikzcd}
\]
Thus $\sigma^*$ is a monomorphism, so \eqref{sturdy complete} holds.

For the reverse direction, it is sufficient to examine pushouts as in \cref{lem pushout stability negative} with $\sigma$ noninvertible and $f$ a generating morphism of $\Tcontr$.
When $f= \sigr$ the condition follows from \eqref{sturdy grounded}, when $f = \delta^\top$ or $\delta^\bot$ the condition follows from \eqref{sturdy stiff}, and when $f=\sigl$ we also have $\sigma = \sigl$ so the condition follows from \eqref{sturdy complete}.
The remaining case is when $f = \sigma^i \colon [n] \twoheadrightarrow [n-1]$, which itself splits into two cases.
If $\sigma(i) = \sigma(i+1)$, then the square factors as below
\[
\begin{tikzcd}
{[n]} \rar[two heads, "\sigma^i"] \dar[swap]{\sigma^i} \drar[phantom, "\ulcorner"very near end] & {[n-1]} \rar[two heads,"\sigma'"] \dar{=}  \drar[phantom, "\ulcorner"very near end] & {[m]} \dar{=}
\\
{[n-1]} \rar[two heads, "="'] & {[n-1]} \rar[two heads, "\sigma'"'] & {[m]}
\end{tikzcd}
\]
and the left hand square is sent to a pullback by \eqref{sturdy complete}.
To handle the remaining case $\sigma(i) \neq \sigma(i+1)$, it suffices to show that 
\[ \begin{tikzcd}
{[m]} \rar{\sigma^k} \dar["\sigma^j"'] \drar[phantom, "\ulcorner" very near end] & {[m-1]} \dar["\sigma^j"] \\
{[m-1]} \rar["\sigma^{k-1}"'] & {[m-2]}
\end{tikzcd} \]
is sent to a pullback whenever $j < k$.
But this now follows from \eqref{sturdy stiff} --
both squares in the following diagram are pushouts 
\[ \begin{tikzcd}
{[1]} \rar{\sigma^0} \dar[tail,"{k,k+1}"'] \drar[phantom, "\ulcorner" very near end] & {[0]} \dar{k}
\\
{[m]} \rar{\sigma^k} \dar["\sigma^j"'] \drar[phantom, "\ulcorner" very near end] & {[m-1]} \dar["\sigma^j"] \\
{[m-1]} \rar["\sigma^{k-1}"'] & {[m-2]}
\end{tikzcd} \]
and the vertical composite $[1] \to [m] \to [m-1]$ on the left is inert.
So \eqref{sturdy stiff} says that the top square and large rectangle are both sent to pullbacks, hence so too is the bottom square.
\end{proof}

\begin{remark}\label{remark stiff complete}
If $X$ is a simplicial space, then condition \eqref{sturdy stiff} for $\lilk^*X$ is called \emph{stiff} in \cite[\S4]{GKT2} (see also \cref{sec:ikeo} below), while condition \eqref{sturdy complete} is called \emph{complete} \cite[\S2.7]{GKT2}.
However, the point of completeness in \cite{GKT2} is to be able to take complements of subobjects, which is not possible when working in more general $\infty$-categories (or even $\infty$-topoi) $\basecat$.
For now we avoid naming condition \eqref{sturdy complete}, and come back to this question in \cref{sec nondeg elts}.
\end{remark}

\section{Ikeo, culf, and conservative maps}
\label{sec map properties}

We extend several types of maps of simplicial objects to presheaves on $\Tsimpl$, $\Tcontr$, and $\Tinert$, taking a cue from the theory of algebraic patterns, a general approach to homotopy coherent algebraic structures \cite{Chu-Haugseng:1907.03977,BarkanHaugsengSteinebrunner:EAP}.
(The definitions we make follow from the notions of active and inert maps in these categories, together with the choice that only $[1]$ is elementary.)

\begin{definition}
A map $Y \to X$ between presheaves over $\Tsimpl$ (or $\Tcontr$ or $\Tinert$) is \emph{ikeo}
if 
\[ \begin{tikzcd}
Y_n \rar \dar \ar[dr,phantom, very near start, "\lrcorner"] & X_n \dar \\
\prod Y_1 \rar & \prod X_1
\end{tikzcd} \]
is a pullback, and $Y_0 \to X_0$ and $Y_\Tob \to X_\Tob$ are equivalences.
\end{definition}

Remembering that the extra object $\Tob$ is 
there to index the zeta function, the ikeo condition for a map of presheaves over $\Tsimpl$ (or $\Tcontr$ or $\Tinert$) thus says ``ikeo in the ordinary sense of simplicial objects'' {\em plus} 
``preservation of zeta functions''.

\begin{remark}\label{rmk rel segal}
In the language of algebraic patterns, ikeo maps would be called ``relative 
$\Tsimpl^{\op,\flat}$-Segal objects'' \cite[Definition 3.1.8]{BarkanHaugsengSteinebrunner:EAP} (where the algebraic pattern $\Tsimpl^{\op,\flat}$ has a single elementary object $[1]$).
\end{remark}

\begin{definition}\label{def generalized culf}
Suppose $\gencopatt$ is a category equipped with an active-inert factorization system.
We say a map $Y \to X$ of $\gencopatt$-presheaves is \emph{culf} if its
restriction to $\gencopatt^\actrm$ is a cartesian natural transformation.
\end{definition}

For the category $\simpl$, this is the usual notion \cite[\S4]{GKT1}.
For $\Tsimpl$, this asks that the underlying map of simplicial objects is 
culf, and that the square
\begin{equation}
\label{eq culf sigr square}
\begin{tikzcd}
Y_\Tob \rar \dar["\sigr^*" swap] & X_\Tob \dar["\sigr^*"] \\
Y_1 \rar & X_1
\end{tikzcd}
\end{equation}
is a pullback (the naturality square for $\sigl^*$ follows from those for
$\sigr^*$ and $s_0$). 
Remembering that the extra object $\Tob$ is 
there to index the zeta function, the culf condition for presheaves over $\Tsimpl$ (or $\Tcontr$ or $\Tinert$) thus says ``culf in the ordinary sense of simplicial objects'' {\em plus} 
``preservation of zeta functions''.

For $\Tcontr$, the culf condition amounts to asking that \eqref{eq culf sigr square} 
is a pullback, as are the squares
\begin{equation}
\label{eq culf degen square}
\begin{tikzcd}
Y_n \rar \dar["s_i" swap] & X_n \dar["s_i"] \\
Y_{n+1} \rar & X_{n+1}.
\end{tikzcd}
\end{equation}
For $\contr$, only the squares \eqref{eq culf degen square} are relevant,
for $\Tinert$ the only square under consideration is \eqref{eq culf sigr square}, 
and for $\inert$ every map is culf.

The pasting law for pullbacks yields the following \cite[Lemma 4.6]{GKT1}.

\begin{lemma}
If $X\in \pre(\Tsimpl)$ is 2-Segal and $Y \to X$ is culf, then $Y$ is 2-Segal. \qed
\end{lemma}

We also have the following weaker condition (see \S\ref{ss ternary}).

\begin{definition}
  A map $Y \to X$ between presheaves over one of the categories in \cref{sec shape
  categories} is \emph{conservative} if it is cartesian when restricted to negative
  maps.
\end{definition}

For $\Tinert$ this is no condition at all: every map in
$\pre(\Tinert)$ is conservative. For $\Tcontr$, a map $Y\to X$ 
in $\pre(\Tcontr)$ is conservative when
\eqref{eq culf degen square} is a pullback and the square
\begin{equation}
\label{eq cons sigl square}
\begin{tikzcd}
Y_0 \rar \dar["\sigl^*" swap] & X_0 \dar["\sigl^*"] \\
Y_\Tob \rar & X_\Tob
\end{tikzcd}
\end{equation}
is a pullback.
For $\Tsimpl$, conservative simply means that the underlying map of
$\Tcontr$-presheaves is conservative, or, in other words, that \eqref{eq cons sigl
square} is a pullback and the underlying map of simplicial objects is conservative
\cite[\S4]{GKT1}. We emphasize that every culf map is conservative, since $\Tsimpl^-
\subset \Tsimpl^\actrm$.

The pasting law for pullbacks gives:
\begin{lemma}\label{conservative sturdy stiff}
If $X$ is sturdy and $Y\to X$ is conservative, then $Y$ is sturdy. \qed
\end{lemma}

\begin{remark}\label{rmk culf vs conservative}
  For $\Tcontr$-presheaves (or $\Tinert$-presheaves), conservative and culf
  maps are not so different: the culf condition is just conservative combined with being
   cartesian on $\sigr \colon [1] \to \Tob$ (as in \eqref{eq culf sigr square}).
  In
  particular, any conservative map between standard $\Tcontr$-presheaves is
  culf.
\end{remark}

\begin{blanko}[The category of \texorpdfstring{$M$}{M}-graded 2-Segal presheaves]\label{M-graded decomp spaces}
  We now spell out the category $\dgr{M}$ of 2-Segal presheaves
  conservatively graded in a $\Tsimpl$-presheaf $M$. This category
  $\dgr{M}$ is precisely the one appearing in \cref{intro theorem B}.  
  While the
  construction of $\dgr{M}$ is general, we eventually wish to show this
  category has a terminal object. In \cref{sec main thm proof} we insist on
  several conditions in order to guarantee this is the case.

Given $M\in \pre(\Tsimpl)$, let $\dcons{M} \subset \pre(\Tsimpl)_{/M}$ denote the full subcategory whose objects are those maps $X\to M$ such that 
\begin{itemize}
\item $X\to M$ a conservative map, and
\item $X$ is 2-Segal.
\end{itemize}
A morphism $Y \to X \to M$ in $\dcons{M}$ will be called ikeo or culf just 
when $Y\to X$ is so.

As ikeo and culf maps are stable under pullback, we have an \emph{adequate triple} 
\begin{equation}\label{eq adequate triple}
  (\dcons{M}, \dcons{M}^\culf, \dcons{M}^\ikeo)
\end{equation}
in the sense of \cite[Definition 5.4]{Barwick:SMFEAK1}.
\end{blanko}

\begin{definition}
The category of $M$-graded 2-Segal presheaves, denoted $\dgr{M}$, is the associated span $\infty$-category of the adequate triple \eqref{eq adequate triple} (see Definition 5.10 of \cite{Barwick:SMFEAK1} or Definition 2.12 of \cite{HaugsengHebestreitLinskensNuiten:2VFFSCS}).
\end{definition}

Concretely, an object of $\dgr{M}$ is a 2-Segal presheaf $X$ equipped with a conservative grading $X \to M$.
Morphisms are commutative diagrams of the form
\[ \begin{tikzcd}
& Z \drar{\culf} \dlar[swap]{\ikeo} & \\
X \drar[swap]{\cons} & & Y \dlar{\cons} \\
& M 
\end{tikzcd} \]
where necessarily $Z$ is again 2-Segal and $Z \to X$ is conservative.
(See also \cite[\S2.1]{Godicke:IC2SS}.)
As is usual with the span $\infty$-category of an adequate triple \cite[3.7, 5.11]{Barwick:SMFEAK1}, the mapping space $\map(f \colon X \to M, g \colon Y \to M)$ is the maximal subgroupoid of the $\infty$-category
\[ 
\dcons{M}_{/f}^\ikeo \times_{\dcons{M}} \dcons{M}_{/g}^\culf \simeq \pre(\Tsimpl)^\ikeo_{/X} \times_{\pre(\Tsimpl)_{/M}} \pre(\Tsimpl)^\culf_{/Y}.
\]
(In the cases we are interested in for \cref{main lemma}, this turns out to be a groupoid.)

\section{Standard objects}\label{section standard}
Taking the clue from standard zeta functions, we call a presheaf $X \in \pre(\Tsimpl)$ \emph{standard} if $\sigr^* \colon X_\Tob \to X_1$ is an equivalence (and likewise for presheaves over $\Tcontr$ or $\Tinert$).
In this section we show that for any object $A \in \pre(\Tcontr)$, the category $\pre(\Tcontr)_{/A}^\std$ (whose objects are those $B\to A$ with $B$ standard) has a terminal object.
This is relatively straightforward.
Slightly more involved is showing that this terminal object is always an ikeo map.
From this we may deduce that $\pre(\Tcontr)_{/A}^{\std,\ikeo}$ is a groupoid, hence contractible.
When $A$ is sturdy, we furthermore show that the terminal object of $\pre(\Tcontr)_{/A}^\std$ is conservative (\cref{lem standarization conservative}).

The two adjunctions
\[
\begin{tikzcd}[column sep=large]
    \pre(\contr)  \rar[hook, "\lilp^*" description] 
    \rar[phantom, bend right=18, "\scriptscriptstyle\perp"]  
    \rar[phantom, bend left=18, "\scriptscriptstyle\perp"] 
    & \pre(\Tcontr)
    \lar[bend right=30, "\lilp_!" swap] 
    \lar[bend left=30, "\lilp_*"] 
\end{tikzcd}
\]
allow us to interpret $\pre(\contr)$ as both a reflective and coreflective subcategory of $\pre(\Tcontr)$. 
As it is coreflective, for each object $A\in \pre(\Tcontr)$ the slice $\pre(\contr)_{/A} \coloneq \pre(\contr) \times_{\pre(\Tcontr)} \pre(\Tcontr)_{/A}$
 has a terminal object given by the counit $\varepsilon_{A} \colon \lilp^* \lilp_* A \to A$ (Proposition 4.1.5 of \cite{RiehlVerity:EICT}).

\begin{lemma}
The image of the fully faithful functor $\lilp^* =\lili_* \colon \pre(\contr) \to \pre(\Tcontr)$ is the subcategory of standard objects. (And similarly for $\lilp^* \colon \pre(\simpl) \to \pre(\Tsimpl)$.)
\end{lemma}
\begin{proof}
$A$ is standard if and only if the unit  $A \to \lilp^*\lilp_! A = \lilp^*\lili^*A$ is an equivalence, and this characterizes the essential image of $\lilp^*$.
\end{proof}

We write $\pre(\Tcontr)^\std_{/A} \subset \pre(\Tcontr)_{/A}$ for the full subcategory whose objects are maps $B\to A$ with standard domain; notice that projection $\pre(\contr) \times_{\pre(\Tcontr)} \pre(\Tcontr)_{/A} \to \pre(\Tcontr)_{/A}$ induces an equivalence $$\pre(\contr)_{/A} \simeq \pre(\Tcontr)^\std_{/A}.$$

We refer to the map $\lilp^* \lilp_* A \to A$ (the terminal object of $\pre(\Tcontr)^\std_{/A}$) as the \emph{standardization
map}.\footnote{Temporary terminology for the process of turning a $\Tsimpl$-presheaf
standard; we regret the clash with the usual terminology ``standardization'' in the
theory of symmetric functions.} Our next aim is to show that it is ikeo. One possible
approach is to rely on the limit formula for pointwise right Kan extension, along with
a cofinality argument. We instead give a slightly more conceptual proof.

Let $\Tgr$ be the full subcategory of $\Tcontr$ on the objects $[0], [1], \Tob$, and $\gr = \Tgr \cap \contr$.
The category $\gr$ is exactly the indexing category for directed reflexive 
graphs, and a presheaf on $\Tgr$ is the same thing as a map of directed 
reflexive graphs that is the identity on vertices. 
Write $\lilu \colon \Tgr \to \Tcontr$ for the inclusion.

The inclusion $\lili \colon \gr \to \Tgr$ has right adjoint $\lilp \colon \Tgr \to \gr$, by restricting the adjunction between the larger categories.
But in our situation, $\lilp$ has a further right adjoint: there is an inclusion
\[ \begin{tikzcd}[row sep=0]
\dirgraph \rar[hook]{\lilw} & \mdirgraph
\end{tikzcd} \]
sending $[0]$ to $[0]$ and $[1]$ to $\Tob$.
We have $\lilp \dashv \lilw$ and the counit of this adjunction is the identity on $\id_{\dirgraph}$. 
Interpreting a $\Tgr$-presheaf as a map of graphs, $\lili^*$ returns the target graph and $\lilw^*$ returns the source graph.
The adjunctions $\lili \dashv \lilp \dashv \lilw$ induce the adjunctions 
\[
\lili_!\quad\dashv \quad \lilp_! = \lili^* \quad\dashv \quad \lilw_! = \lilp^* = \lili_* \quad\dashv \quad \lilw^* = \lilp_* \quad\dashv \quad \lilw_*
\]
between $\pre(\gr)$ and $\pre(\Tgr)$, with the outer and middle functors fully faithful.
The primary benefit to the current approach is that $\lilp_*$ at this level becomes entirely transparent.

\begin{proposition}
Consider the adjunction $\lilu^* : \pre(\Tcontr) \rightleftarrows \pre(\Tgr) : \lilu_*$ and let $A \to \bar A \coloneq \lilu_*\lilu^* A$ be the unit of the adjunction at $A$.
The standardization map $\lilp^*\lilp_* A \to A$ is given by the projection $(\lilp^*\lilp_*\bar A) \times_{\bar A} A \to A$.
\end{proposition}
\begin{proof}
As $\lilp^* \lilp_*$ preserves limits and is idempotent, the pullback below left is sent to the pullback below right.
\[ \begin{tikzcd}
P \rar \drpullback \dar & A \dar \\
\lilp^* \lilp_* \bar A \rar & \bar A
\end{tikzcd} 
\rightsquigarrow
\begin{tikzcd}
\lilp^*\lilp_* P \rar \drpullback \dar & \lilp^*\lilp_* A \dar \\
\lilp^* \lilp_* \bar A \rar{\simeq} & \lilp^*\lilp_* \bar A
\end{tikzcd} 
\]
Since the bottom map to the right is an equivalence, so too is the top map, and it suffices to show that $P$ is standard.
But $\lilu^* (A \to \bar A)$ is an equivalence, hence so is $\lilu^*(P \to \lilp^*\lilp_* \bar A)$.
Thus $P$ is standard.
\end{proof}

Now $\lilp^* \lilp_* \bar A = \lilp^* \lilp_* \lilu_* \lilu^* A$. The adjunction
formulas above give $\lilp^* \lilp_* \lilu_* = \lili_* \lilp_* \lilu_* = \lilu_*
\lili_* \lilp_* = \lilu_* \lilp^* \lilp_*$ and the adjunction counit $\lilp^*
\lilp_* \bar A \to \bar A$ is the image of the adjunction counit $\lilp^* \lilp_*
\lilu^* A \to \lilu^*A$ under $\lilu_*$.

\begin{lemma}\label{lem warmup ikeo}
If $V' \to V$ is a map in $\pre(\Tgr)$ that is an equivalence at $[0]$ and $\Tob$, then $\lilu_* V' \to \lilu_* V$ is ikeo.
\end{lemma}
\begin{proof}
A calculation gives $(\lilu_*V)_n \simeq V_1 \times_{V_0} \cdots \times_{V_0} V_1$, so $\lili^* \lilu_*V \in \pre(\contr)$ is Segal.
Since $V'_0 \to V_0$ is an equivalence, $\lili^* \lilu_* (V' \to V)$ is ikeo.
If $f \colon V' \to V$ is an equivalence at $[0]$ and $\Tob$, then so is $\lilu_* V' \to \lilu_* V$; on underlying $\contr$-presheaves this map is then ikeo since both source and target are Segal.
\end{proof}

\begin{proposition}\label{prop standardization ikeo}
The standardization map $\lilp^*\lilp_* A \to A$ is ikeo.
\end{proposition}
\begin{proof}
Let $V = \lilu^* A$. 
The map $\lilp^* \lilp_* V \to V$ is $(\lilw \lilp)^* V \to V$ coming from the unit $\id \Rightarrow \lilw \lilp$, which is the identity at $[0]$ and $\Tob$.
Thus $\lilp^* \lilp_* V \to V$ is the identity at $[0]$ and $\Tob$.
Applying $\lilu_*$ to this gives an ikeo map by \cref{lem warmup ikeo}.
\end{proof}

\begin{lemma}\label{lem standarization conservative}
If $A\in \pre(\Tcontr)$ is grounded, then $\lilp^*\lilp_* A \to A$ is conservative.
\end{lemma}
\begin{proof}
Since conservative maps are stable under pullback, it suffices to show that $\lilp^* \lilp_* \bar A \to \bar A$ is conservative.
For this, it suffices by \cref{lem warmup ikeo} to show that $\lilp^* \lilp_* \lilu^* A \to \lilu^* A$ is conservative in $\pre(\Tgr)$, namely to check that the squares 
\[ \begin{tikzcd}
\lilp^* \lilw^* \lilu^* A_0 \rar \dar & \lilp^* \lilw^* \lilu^* A_1 \dar  
& 
\lilp^* \lilw^* \lilu^* A_0 \rar \dar & \lilp^* \lilw^* \lilu^* A_\Tob \dar  \\
\lilu^* A_0 \rar & \lilu^* A_1 
&
\lilu^* A_0 \rar & \lilu^* A_\Tob
\end{tikzcd} \]
are pullbacks.
But the left square is a pullback since $A$ is grounded (\cref{prop sturdy char}(\ref{sturdy grounded})), and in the right square both vertical maps are equivalences.
\end{proof}

\begin{lemma}\label{lem ikeo standard}
Any ikeo map between standard objects is an equivalence.
\end{lemma}
\begin{proof}
Suppose $D\to E$ is such a map; by ikeo, $D_\Tob \to E_\Tob$ is an  equivalence. 
By standard, $D_\Tob \to D_1$ and $E_\Tob \to E_1$ are equivalences, hence $D_1 \to E_1$ is an equivalence.
Any ikeo map which is an equivalence on 1-simplices is itself an equivalence.
\end{proof}

Since ikeo maps are closed under left cancellation, we deduce the following.

\begin{proposition}\label{prop std ikeo contractible}
If $A$ is a grounded $\Tcontr$-presheaf (e.g.\ if $A$ is sturdy), then the category $\pre(\Tcontr)_{/A}^{\std,\ikeo,\cons}$ is a contractible groupoid.
\end{proposition}
\begin{proof}
  The category $\pre(\Tcontr)_{/A}^{\std,\ikeo,\cons}$ is a full subcategory of
  $\pre(\Tcontr)_{/A}^{\std,\ikeo}$ which is a groupoid by \cref{lem ikeo standard}. 
  The terminal object of $\pre(\Tcontr)_{/A}^{\std}$ 
  is a map which is
  both ikeo and conservative by \cref{prop standardization ikeo} and \cref{lem
  standarization conservative}, so it belongs to 
  $\pre(\Tcontr)_{/A}^{\std,\ikeo,\cons}$,
  which is therefore a groupoid with a terminal object, hence contractible. 
\end{proof}

\section{Left Kan extension along \texorpdfstring{$\lilh$}{h}, \texorpdfstring{$\lilj$}{j}, and \texorpdfstring{$\lilk$}{k}}\label{sec lke along jkl}

In this section we explore left Kan extension along the functors $\lilh$, $\lilj$, and $\lilk$:
\[ \begin{tikzcd}[column sep=0, row sep=small]
\Tinert \ar[rr,"\lilj"] \drar["\lilh"'] & & \Tsimpl \\
& \Tcontr \urar["\lilk"']
\end{tikzcd} \]
We have much more to say here about the first two since these admit nice formulas.
Studying $\lilk_!$ is more difficult, and partial results in \cref{sec culf conserv lke} and \cref{sec:ikeo} rely on our study of $\lilh_!$ and $\lilj_!$ here.

The three categories $\Tinert$, $\Tcontr$, and $\Tsimpl$ contain exactly the same subcategory of inert maps $\Tinert^\intrm = \Tcontr^\intrm = \Tsimpl^\intrm = \inert \cup \{\id_\Tob, \sigr \delta^0, \sigr \delta^1 \}$.
This allows us to reduce the usual formulas for left Kan extension to only consider active maps.
Let us examine how this works for $\lilj\colon \Tinert \to \Tsimpl$.
For $e\in \Tsimpl$, the map on coslices
\[
\Tsimpl^\actrm_{e/} \to \Tsimpl_{e/}
\]
has a right adjoint given by active-inert factorization.
Since $\Tinert$ and $\Tsimpl$ have the same class of inert maps, this induces an adjunction
\[
        \Tinert^\actrm \times_{\Tsimpl^\actrm} \Tsimpl^\actrm_{e/} \rightleftarrows  \Tinert \times_{\Tsimpl} \Tsimpl_{e/}.
\]
(The objects in these comma categories are (active) maps $e \ract d$ in $\Tsimpl$ from the fixed object $e\in\Tsimpl$ to some arbitrary object in $\Tinert$.)
In particular, the inclusion $\Tinert^\actrm \times_{\Tsimpl^\actrm} \Tsimpl^\actrm_{e/} \to  \Tinert \times_{\Tsimpl} \Tsimpl_{e/}$ is an initial functor, so
\[
        (\lilj_! A)_e = \colim_{e \ract d} A_d
\]
where the colimit is taken over the category $(\Tinert^\actrm \times_{\Tsimpl^\actrm} \Tsimpl^\actrm_{e/})^\op$.
(One could instead invoke \cite[Lemma 7.2]{Chu-Haugseng:1907.03977}.)
Similar considerations hold for computing the left Kan extension along $\lilh\colon \Tinert \to \Tcontr$ or $\lilk \colon \Tcontr \to \Tsimpl$.

\begin{lemma}[Beck--Chevalley]\label{beck-chevalley}
The natural transformations
\begin{align*}
\lilk_! \lili^* &\Rightarrow \lili^* \lilk_!  \colon \pre(\Tcontr) \to \pre(\simpl) \\
\lilh_! \lili^* &\Rightarrow \lili^* \lilh_!  \colon \pre(\Tinert) \to \pre(\contr) \\
\lilj_! \lili^* &\Rightarrow \lili^* \lilj_!  \colon \pre(\Tinert) \to \pre(\simpl)
\end{align*}
are equivalences.
\end{lemma}
\begin{proof}
The inclusion $\lili$ has a retraction $\lilp$ with $\lili \dashv \lilp$, so we have $\lili^* \simeq \lilp_!$.
The results now follow using commutativity of \eqref{commutative ipjlk} just below \cref{notation for functors}.
\end{proof}

\begin{lemma}\label{lem k shriek at T}
If $A\in \pre(\Tcontr)$, then $(\lilk_!A)_\Tob = A_\Tob$.
\end{lemma}
\begin{proof}
This follows immediately from the fact that $\id_\Tob$ is the initial object of 
\[
        \Tcontr^\actrm \times_{\Tsimpl^\actrm} \Tsimpl^\actrm_{\Tob/}.
\]
There are only two objects in this category, since the only other active map of $\Tsimpl$ with domain $\Tob$ is $\sigl \colon \Tob \to [0]$.
But $\sigl$ is in $\Tcontr$, so indeed $\id_\Tob$ is an initial object.
\end{proof}

\begin{proposition}\label{lke of hj}
If $A\in \pre(\Tinert)$ then 
\begin{enumerate}
\item $(\lilh_! A)_\Tob =(\lilj_! A)_\Tob = A_0 + A_\Tob$ \label{item hj lke at T}
\item  \label{item hj lke away from T} For an object $[n]$, we have \[(\lilh_! A)_n = \sum_{[n] \twoheadrightarrow [k]} A_k \qquad \text{ and } \qquad  (\lilj_! A)_n = \sum_{[n] \ract [k]} A_k.\]
\end{enumerate}
\end{proposition}
\begin{proof}
For \eqref{item hj lke at T}, the categories $\Tinert^\actrm \times_{\Tcontr^\actrm} \Tcontr^\actrm_{\Tob/}$ and $\Tinert^\actrm \times_{\Tsimpl^\actrm} \Tsimpl^\actrm_{\Tob/}$ are discrete categories on two objects $\id_\Tob$ and $\sigl \colon \Tob \to [0]$.
Discreteness is because there are no maps in $\Tinert^\actrm$ between $\Tob$ and $[0]$.

By \cref{beck-chevalley}, to establish \eqref{item hj lke away from T}, it suffices to compute $\lilh_!$ and $\lilj_!$ on $\lili^*A \in \pre(\inert)$.
We may compute these left Kan extensions using the categories $\inert^\actrm \times_{\contr^\actrm} \contr^\actrm_{[n]/} = \inert^\actrm \times_{\simpl^-} \simpl^-_{[n]/}$ and $\inert^\actrm \times_{\simpl^\actrm} \simpl^\actrm_{[n]/}$, respectively.
But since $\inert^\actrm$ is discrete, so too are these categories: $\inert^\actrm \times_{\contr^\actrm} \contr^\actrm_{[n]/} = \{[n] \twoheadrightarrow [k]\}$ and  $\inert^\actrm \times_{\simpl^\actrm} \simpl^\actrm_{[n]/} = \{[n] \ract [k]\}$.
\end{proof}
Notice that \eqref{item hj lke away from T}  mimics \cite[Corollary 2.1.4]{Hackney-Kock:2210.11192}.

\begin{remark}\label{rem pullback preservation}
The preceding proposition implies that $\lilh_!$ and $\lilj_!$ preserve pullbacks, since pullbacks commute with sums in $\basecat$.
\end{remark}

\section{Culf and conservative maps under left Kan extension}\label{sec culf conserv lke}

Given objects $A\in \pre(\Tinert)$, $B\in \pre(\Tcontr)$, we have induced left adjoints on slice categories $\lilh_! \colon \pre(\Tinert)_{/A} \to \pre(\Tcontr)_{/\lilh_!A}$ and $\lilk_! \colon \pre(\Tcontr)_{/B} \to \pre(\Tsimpl)_{/\lilk_!B}$ and $\lilj_! \colon \pre(\Tinert)_{/A} \to \pre(\Tsimpl)_{/\lilj_!A}$.
One might wonder whether these restrict to the full subcategories whose objects are culf or conservative maps.
Our goal for this section is to prove the following theorem, which is prefigured in \cite[Theorem 3.2.1]{Hackney-Kock:2210.11192}.

\begin{theorem}\label{thm four equivs}
  If $A \in \pre(\Tinert)$, then we have the following induced functors
  \[ \begin{tikzcd}
  \pre(\Tinert)^\culf_{/A} \rar{\lilh_!} \ar[rr, bend right=15, "\lilj_!"'] 
  & \pre(\Tcontr)^\culf_{/\lilh_!A} \rar{\lilk_!} & \pre(\Tord)^\culf_{/\lilk_!\lilh_!A}.
  \end{tikzcd} \]
  These are all equivalences.
  Likewise, there is an equivalence $\lilh_! \colon \pre(\Tinert)_{/A} 
  = \pre(\Tinert)^\cons_{/A} \simeq \pre(\Tcontr)^\cons_{/\lilh_!A}$.
\end{theorem}
To prove the statement about $\lilk_!$, it suffices to work entirely with $\lilh_!$ and $\lilj_!$.
We do not know whether one should expect to have an equivalence $\lilk_! \colon \pre(\Tcontr)^\culf_{/B} \simeq \pre(\Tord)^\culf_{/\lilk_!B}$ for general $B$.

Many formulas below involve pullbacks of sums over sets, and we implicitly use the following several times.
\begin{lemma}\label{lem little coproducts lemma}
Suppose $V$ is a set and we have $V$-indexed collections of objects $A_v, B_v$ in $\basecat$, along with a function $\psi \colon U \to V$.
Given a $V$-indexed collection of maps $B_v \to A_v$, the associated square
\[ \begin{tikzcd}
\displaystyle \sum_{u\in U} B_{\psi(u)} \rar \dar & \displaystyle \sum_{v\in V} B_v \dar \\
\displaystyle \sum_{u\in U} A_{\psi(u)} \rar & \displaystyle \sum_{v\in V} A_v 
\end{tikzcd} \]
is a pullback.
\end{lemma}
\begin{proof}
Write $V'$ for the object $\sum_{v\in V} 1$, so a $V$-indexed collection is a map $B \to V'$ (and similarly for $U'$).
Then the square above is the top square below
\[ \begin{tikzcd}
\psi^* B \rar \dar & B \dar \\
\psi^* A \rar \dar \drar[phantom, "\lrcorner" very near start] & A \dar \\
U' \rar{\psi} & V',
\end{tikzcd} \]
which is thus a pullback by the pasting law.
\end{proof}

\begin{lemma}\label{lem culf preservation}
The functors $\lilh_!$ and $\lilj_!$ take values in conservative maps.
A map $f$ in $\pre(\Tinert)$ is culf if and only if $\lilh_!(f)$ is culf if and only if $\lilj_!(f)$ is culf.
\end{lemma}
\begin{proof}
Suppose $f\colon B\to A$ is a map between $\Tinert$-presheaves, and $\sigma \colon [n] \twoheadrightarrow [m]$ a map in $\simpl^- = \contr^\actrm$.
We wish to show that $\lilh_!(f)$ is cartesian at $\sigma$.
By \cref{lke of hj}\eqref{item hj lke away from T} the squares below are the same 
\[ 
\begin{tikzcd}
\lilh_! B_m \rar \dar & \lilh_! B_n \dar \\
\lilh_! A_m \rar & \lilh_! A_n
\end{tikzcd} 
\qquad \text{and} \qquad
\begin{tikzcd}
\displaystyle \sum_{[m] \twoheadrightarrow [k]} B_k \rar \dar & \displaystyle \sum_{[n] \twoheadrightarrow [k]} B_k \dar \\
\displaystyle \sum_{[m] \twoheadrightarrow [k]} A_k \rar & \displaystyle \sum_{[n] 
\twoheadrightarrow [k]} A_k .
\end{tikzcd} 
\]
But the right square is a pullback by \cref{lem little coproducts lemma}.
Similarly, $\lilj_!(f)$ is cartesian at maps of $\simpl^\actrm$  (or see \cite[Proposition 2.3.2]{Hackney-Kock:2210.11192}), so in particular at maps of $\simpl^-$.
There is one more map to worry about, namely $\sigl \colon \Tob \to [0]$.
But by \cref{lke of hj}\eqref{item hj lke at T}, the maps $\lilh_!(f)$ or $\lilj_!(f)$ at $\sigl$ are both the square
\[ \begin{tikzcd}
B_0 \rar \dar & B_0 + B_\Tob \dar \\
A_0 \rar & A_0 + A_\Tob
\end{tikzcd} \]
which is a pullback.
Hence $\lilh_!(f)$ and $\lilj_!(f)$ are conservative.

The map $f\colon B \to A$ is culf if and only if it is cartesian at $\sigr \colon [1] \to \Tob$.
The naturality square for $\lilj_!(f)$ at $\sigr$ is the outer rectangle below (by \cref{lke of hj}):
\begin{equation} \label{eq sigr jshriek}
\begin{tikzcd}
B_0 + B_\Tob \rar{\id + \sigr^*} \dar & B_0 + B_1 \rar[hook] \dar \drpullback & \displaystyle \sum_{k\in \mathbb{N}} B_k \dar \\
A_0 + A_\Tob \rar{\id + \sigr^*} & A_0 + A_1 \rar[hook] & \displaystyle \sum_{k\in \mathbb{N}} A_k.
\end{tikzcd} \end{equation}
Thus the left square is a pullback if and only the outer rectangle is a pullback. 
It follows that $f$ is cartesian at $\sigr$ if and only if $\lilj_!(f)$ is so.
Since $\lilj_!(f)$ is automatically cartesian at all active maps in $\simpl$, this establishes that $f$ is culf if and only if $\lilj_!(f)$ is culf.
The argument for $\lilh_!(f)$ is similar, just noting that the left square of \eqref{eq sigr jshriek} is exactly the naturality square for $\lilh_!(f)$ at $\sigr$.
\end{proof}

\begin{lemma}\label{lem cartesian unit}
The units for the $\lilh_! \dashv \lilh^*$ and $\lilj_! \dashv \lilj^*$ adjunctions are cartesian.
\end{lemma}
\begin{proof}
This essentially follows from \cref{lke of hj}.
For the $\lilh$ case, the unit $\eta \colon \id \Rightarrow \lilh^* \lilh_!$ is cartesian since the following squares are pullbacks for any map $B\to A$ of $\Tinert$-presheaves. 
\[ \begin{tikzcd}
B_n \rar[hook] \dar \drpullback & \displaystyle  \sum_{[n] \twoheadrightarrow [k]} B_k \dar & B_\Tob \rar[hook] \dar \drpullback & B_0 + B_\Tob \dar \\
A_n \rar[hook] & \displaystyle  \sum_{[n] \twoheadrightarrow [k]} A_k & A_\Tob \rar[hook] & A_0 + A_\Tob
\end{tikzcd} 
\]
A similar argument establishes that $\id \Rightarrow \lilj^*\lilj_!$ is cartesian.
\end{proof}

\begin{lemma}\label{lem counit cartesian}
  The counit for the $\lilh_! \dashv \lilh^*$ adjunction is cartesian at conservative
  maps in $\pre(\Tcontr)$ (and hence at culf maps) and the counit for $\lilj_! \dashv \lilj^*$ is
  cartesian at culf maps in $\pre(\Tsimpl)$.
\end{lemma}
\begin{proof}

Suppose $B \to A$ is a conservative map of $\Tcontr$-presheaves.
Then for each degeneracy $\sigma \colon [n] \twoheadrightarrow [k]$, the square below left is a pullback.
\[ \begin{tikzcd}
B_k \rar{\sigma^*}\dar \drpullback & B_n \dar & \displaystyle  \sum_{[n] \twoheadrightarrow [k]} B_k \rar\dar & B_n \dar\\
A_k \rar[swap]{\sigma^*} & A_n & \displaystyle  \sum_{[n] \twoheadrightarrow [k]} A_k \rar & A_n
\end{tikzcd} \]
It follows that the square above right is a pullback, but this is precisely the naturality square for the counit at $n$.
At $\Tob$, the naturality square is (by the proof of \cref{lke of hj})
\[ \begin{tikzcd}
B_0 + B_\Tob \rar{\sigl^* + \id} \dar & B_\Tob \dar \\
A_0 + A_\Tob \rar[swap]{\sigl^* + \id} & A_\Tob ,
\end{tikzcd} \] 
which is cartesian since $B \to A$ is cartesian at $\sigl \in \Tcontr^-$.
The case for $\varepsilon \colon \lilj_! \lilj^* \Rightarrow \id$ is analogous.
\end{proof}

Recall that given an adjunction of $\infty$-categories $F : \mathcal{C} \rightleftarrows \mathcal{D} : G$ where $\mathcal{C}$ admits pullbacks, then for every object $c\in \mathcal{C}$ there is an induced adjunction $\mathcal{C}_{/c} \rightleftarrows \mathcal{D}_{/Fc}$ whose right adjoint is given by applying $G$ and pulling back along the unit natural transformation: $G(d\to Fc) = (c \times_{GFc} Gd \to c)$ \cite[Proposition~5.2.5.1]{Lurie:HTT}.
A variation of the proof of \cite[\S8.9]{Hackney-Kock:2210.11191} gives the following more precise version of \cref{thm four equivs}.

\begin{theorem}\label{thm restricted slices}
The sliced adjunctions 
\begin{xalignat*}{2}
\pre(\Tinert)_{/A} &\rightleftarrows \pre(\Tcontr)_{/\lilh_!A} &
 \pre(\Tinert)_{/A} &\rightleftarrows \pre(\Tsimpl)_{/\lilj_!A}
\intertext{
restrict to adjoint equivalences 
}
\pre(\Tinert)^\culf_{/A} &\rightleftarrows \pre(\Tcontr)^\culf_{/\lilh_!A} 
& \pre(\Tinert)^\culf_{/A} &\rightleftarrows \pre(\Tsimpl)^\culf_{/\lilj_!A} \\
 \pre(\Tinert)_{/A} &\rightleftarrows \pre(\Tcontr)^\cons_{/\lilh_!A}. &&
\end{xalignat*}
\end{theorem}
\begin{proof}
We focus on the first equivalence, as the other two use the same formal argument.
The left adjoint in the sliced adjunction $\pre(\Tinert)_{/A} \rightleftarrows \pre(\Tcontr)_{/\lilh_!A}$ is just given by $\lilh_!$, while the right adjoint is given by first restricting with $\lilh^*$ and then pulling back along $\eta_A \colon A \to \lilh^* \lilh_!A$.
Since $\lilh^*$ clearly preserves culf maps and $\eta_A$ is cartesian (\cref{lem cartesian unit}), the right adjoint sends culf maps (considered as objects in the slice) to culf maps.
Using \cref{lem culf preservation}, we thus have that the sliced adjunction restricts to the full subcategories of culf maps: $\pre(\Tinert)^\culf_{/A} \rightleftarrows \pre(\Tcontr)^\culf_{/\lilh_!A}$.

We show that the counit and unit of this adjunction are equivalences.
We start with the counit.
For a culf map $X \to \lilh_!A$ between $\Tcontr$-presheaves, the value of the right adjoint is the left map in the pullback square
\[ \begin{tikzcd}
B \dar \rar \drpullback & \lilh^* X \dar \\
A \rar[swap]{\eta_A} & \lilh^*\lilh_! A  .
\end{tikzcd} \]
The counit is given by the following map over $\lilh_!A$:
\[ \begin{tikzcd}
\lilh_! B \dar \rar \drpullback & \lilh_! \lilh^* X \rar{\varepsilon_X} \dar \drpullback & X \dar \\ 
\lilh_! A \rar{\lilh_!\eta_A} \ar[rr, bend right=15,"\id"'] & \lilh_! \lilh^* \lilh_! A \rar{\varepsilon_{\lilh_! A}} & \lilh_! A.
\end{tikzcd} \]
The left square is a pullback by \cref{rem pullback preservation}, and the right square is a pullback by \cref{lem counit cartesian}, hence $\lilh_!B \to X$ is an equivalence.
Turning to the unit of the adjunction, we simply note that for any map $B\to A$, the square
\[ \begin{tikzcd}
B \rar{\eta_B} \dar \drpullback & \lilh^* \lilh_! B \dar \\
A \rar[swap]{\eta_A} & \lilh^* \lilh_! A
\end{tikzcd} \]
is a pullback by \cref{lem cartesian unit}.
\end{proof}

\section{The \texorpdfstring{$\Tinert$}{Ĩ}-presheaf of nondegenerate elements}\label{sec nondeg elts}

In this section we study \emph{completely sturdy} $\Tcontr$-presheaves, namely those of the form $\lilh_!\nondeg A$ for an $\Tinert$-presheaf $\nondeg A$ of nondegenerate elements (\cref{def nondeg A}).
These are exactly the sturdy presheaves whose degeneracy maps are complementable:

\begin{proposition}\label{prop three viewpoints}
For $A \in \pre(\Tcontr)$, the following are equivalent:
\begin{enumerate}
\item $A$ is in the essential image of $\lilh_! \colon \pre(\Tinert) \to \pre(\Tcontr)$, \label{item:3vp-image}
\item $A$ admits a conservative map to $\lilh_! 1$, and \label{item:3vp-cons}
\item $A$ is sturdy and all degeneracies (including $\sigl^*$) are sum inclusions.\label{item:3vp-degen}
\end{enumerate}
\end{proposition}

\begin{definition}\label{def completely sturdy}
We say that $A \in \pre(\Tcontr)$ is \emph{completely sturdy} when it satisfies the equivalent conditions from \cref{prop three viewpoints}.
We also say $X\in \pre(\Tsimpl)$ is completely sturdy when $\lilk^*X$ is so.
\end{definition}

Notice that when $\basecat = \spaces$, every sturdy presheaf is completely sturdy.

\begin{proof}[Proof of \cref{prop three viewpoints}]
By \cref{thm restricted slices}, the functor $\lilh_!$ induces an equivalence 
\[
   \pre(\Tinert) \simeq \pre(\Tinert)_{/1} \to \pre(\Tcontr)^\cons_{/\lilh_!1},
\]
yielding the equivalence between \eqref{item:3vp-image} and \eqref{item:3vp-cons}.
Meanwhile, \eqref{item:3vp-cons} will be shown to be equivalent to \eqref{item:3vp-degen} in \cref{lem degeneracy characterization} below.
\end{proof}

Our first task will be to show that \eqref{item:3vp-cons} can be strengthened to ask for a \emph{unique} conservative map to $\lilh_!1$. 
We observe from \cref{lke of hj} that $(\lilh_! 1)_d = \sum_{d\twoheadrightarrow c} 1_c$.

\begin{lemma}\label{lem at most one conservative}
If $A\in \pre(\Tcontr)$, there is at most one conservative map $A \to \lilh_! 1$.
\end{lemma}
\begin{proof}
In this proof we use that, by extensivity of the base $\basecat$, sum inclusions are stable under pullback and complements of subobjects are unique when they exist.
Given any $f\colon A \to \lilh_! 1$, the map $A_d\to (\lilh_! 1)_d$ splits $A_d$ into a sum $\sum_{d\twoheadrightarrow c} A_\sigma^f$ 
where $A_\sigma^f \subset A_d$ is given by the pullback
\[ \begin{tikzcd}
A_\sigma^f \rar[hook] \dar \drar[phantom, "\lrcorner" very near start]& A_d \dar{f} \\
1_{c'} \rar[hook,"\name{\sigma}"] & \displaystyle\sum_{d\twoheadrightarrow c} 1_c
\end{tikzcd} \]
for each $\sigma \colon d \twoheadrightarrow c'$ (this is just an instance of the equivalence $\basecat_{/\sum_{u\in U} 1} \simeq \prod_{u\in U} \basecat$).
Assume now that $f$ is conservative. For any $\tau \colon d' \twoheadrightarrow d$ the degeneracy $A_d \to A_{d'}$ restricts to an equivalence $A_\sigma^f \simeq A_{\sigma\tau}^f$,
since both these objects occupy the $\cdot$ position in the pullback diagram
\[ \begin{tikzcd}
\cdot \rar[hook] \dar \drar[phantom, "\lrcorner" very near start] & A_d \rar{\tau^*} \dar \drar[phantom, "\lrcorner" very near start] & A_{d'} \dar  \\
1_{c'} \rar{\name{\sigma}} \ar[rr,bend right=15,"\name{\sigma\tau}"' near start] & 
(\lilh_! 1)_d \rar{\tau^*} & (\lilh_! 1)_{d'} .
\end{tikzcd} \]

Given another conservative map $g \colon A \to \lilh_! 1$, we will show that $A_\sigma^f = A_\sigma^g$ inside $A_d$ (technically, by this equality we mean that the monomorphisms $A_\sigma^f \to A_d$ and $A_\sigma^g \to A_d$ are equivalent in $\basecat_{/A_d}$, and similarly below).
We will then have that $f$ and $g$ agree at the object $d$ since $f, g\colon A_d \to (\lilh_! 1)_d$ preserve the sum splitting.
As $(\lilh_! 1)_0$ is a point, $A_{\id_{[0]}}^f = A_0 = A_{\id_{[0]}}^g$.
By the previous paragraph, we have $A^f_\sigma = A^g_\sigma$ inside of $A_d$ for all $\sigma\colon d \twoheadrightarrow [0]$.
We then have 
\[
  A_0 + A_{\id_\Tob}^f = A_{\sigl}^f + A_{\id_\Tob}^f = A_\Tob = A_{\sigl}^g + A_{\id_\Tob}^g =  A_0 + A_{\id_\Tob}^g,
\]
so $A_{\id_\Tob}^f = A_{\id_\Tob}^g$.

To see that $A^f_\sigma = A^g_\sigma$ inside of $A_n$ for $\sigma \colon [n] \twoheadrightarrow [k]$, we induct on $k$.
The base case $k=0$ is above.
Inductively, suppose this is true for all $k < n$, 
that is, 
$A_\sigma^f = A_\sigma^g$ inside $A_n$ for each noninvertible surjection $\sigma \colon [n] \twoheadrightarrow [k]$.
But then 
\[
  A_{\id_{[n]}}^f + \sum_{\substack{\sigma \colon [n] \twoheadrightarrow [k]\\\text{noninvertible}}} A_\sigma^f = A_n = A_{\id_{[n]}}^g + \sum_{\substack{\sigma \colon [n] \twoheadrightarrow [k]\\\text{noninvertible}}} A_\sigma^g
\]
and we conclude $A_{\id_{[n]}}^f = A_{\id_{[n]}}^g$ since complements are unique.
By the first paragraph, $A_\tau^f = A_\tau^g$ inside $A_m$ for all $\tau \colon [m] \twoheadrightarrow [n]$.
\end{proof}

\begin{definition}\label{def nondeg A}
If $A \in \pre(\Tcontr)$ admits a (necessarily unique) conservative map $f$ to $\lilh_!1$ (i.e.\ $A$ is completely sturdy), we define $\nondeg A \in \pre(\Tinert)$ to be the unique presheaf such that $\lilh_!(\vec{A} \to 1) = f$.
Explicitly, by \cref{thm restricted slices}, $\vec{A}$ sits in the following pullback square
\[ \begin{tikzcd}
\nondeg A \rar{\eta_{\nondeg A}} \dar \drar[phantom,"\lrcorner" very near start] & \lilh^* A \simeq \lilh^* \lilh_! \nondeg A \dar{\lilh^* f} \\
1 \rar{\eta_{1}} & \lilh^* \lilh_! 1 .
\end{tikzcd} \]
If $X\in \pre(\Tsimpl)$ is completely sturdy with $A \coloneq \lilk^*X$, then we write $\nondeg X \coloneq \nondeg A$.
\end{definition}

\begin{remark}
We have a sum splitting
\[
  A_d = \lilh_!(\vec{A})_d =  \sum_{\sigma \colon d \twoheadrightarrow c} \vec{A}_c = \vec{A}_d + \mkern-25mu \sum_{\substack{\sigma \colon d \twoheadrightarrow c\\\text{noninvertible}}} \vec{A}_c;
\]
the subobject $\vec{A}_d \subset A_d$ of nondegenerate simplices was called $A^f_{\id_{d}}$ in the proof.
One should interpret the formula above as an Eilenberg--Zilber lemma.
In case $d = \Tob$, this formula gives a sum splitting $A_\Tob = A_0 + \vec{A}_\Tob$ whose summands are, respectively, the pullbacks
\[
\begin{tikzcd}[row sep={4em,between origins}, column sep={5em,between origins}]
A_0 \ar[d] \ar[r] \drpullback  & A_\Tob \ar[d]  \\
A_0 \ar[r] & A_1
\end{tikzcd}
\qquad
\begin{tikzcd}[row sep={4em,between origins}, column sep={5em,between origins}]
\nondeg A_\Tob \ar[d] \ar[r] \drpullback  & A_\Tob \ar[d]  \\
\nondeg A_1 \ar[r] & A_1 .
\end{tikzcd}
\]
\end{remark}

\begin{lemma}\label{lilh 1 CSS}
$\lilh_! 1$ is sturdy, and all its degeneracies are sum inclusions.\qed
\end{lemma}

\begin{proposition}\label{lem degeneracy characterization}
A presheaf $A \in \pre(\Tcontr)$ admits a conservative map to $\lilh_! 1$ if and only if $A$ is sturdy and all degeneracies (including $\sigl^*$) are sum inclusions.
\end{proposition}
\begin{proof}
Suppose $A \to \lilh_! 1$ is conservative. Since $\lilh_! 1$ is sturdy by 
\Cref{lilh 1 CSS},
it follows from \cref{conservative sturdy stiff} that also $A$ is sturdy.
Since the degeneracies of $\lilh_!1$ are sum inclusions,
and since sum inclusions are stable under pullback in $\basecat$, it follows that also the degeneracies in $A$ are sum inclusions.

Now suppose $A$ is sturdy and all degeneracies are sum inclusions. 
One may imitate the $1\Rightarrow 2$ direction of \cite[Proposition 4.7]{GKT2} to produce an associated $\pre(\inert)$-presheaf, which is part of a $\pre(\Tinert)$-presheaf we call $\vec{A}^\textup{alt}$ (soon to be shown equivalent to $\nondeg A$).
We define $\vec{A}_\Tob^\textup{alt} \subset A_\Tob$ to be the complement of the subobject $\sigl^* \colon A_0 \to A_\Tob$.
Notice that the restriction of $\sigr^*$ to $\vec{A}_\Tob^\textup{alt}$ lands in the nondegenerate part $\vec{A}_1^\textup{alt} \subset A_1$ since $A$ is grounded. 
Altogether we have constructed a subpresheaf $\vec{A}^\textup{alt} \subset \lilh^*A$ in $\pre(\Tinert)$ consisting of the nondegenerate simplices.
That the adjunct $\lilh_! (\vec{A}^\textup{alt}) \to A$ is an equivalence is a consequence of 
\cref{lke of hj} and a variation of \cite[Corollary 2.17]{GKT2} (so in particular, $\vec{A}^\textup{alt} \simeq \vec{A}$ by \cref{thm restricted slices}). 
Since $A$ is thus in the image of $\lilh_!$, it admits a conservative map 
 to $\lilh_!1$.  
\end{proof}

\begin{remark}\label{rmk functoriality of nondeg}
The assignment $A \mapsto \nondeg A$ is a functor (in fact an equivalence) from the category of completely sturdy $\Tcontr$-presheaves and conservative maps to $\pre(\Tinert)$.
This follows from \cref{thm restricted slices}, the definition of $\nondeg A$, and the pasting law for pullbacks.
If $A \to B$ is a conservative map between completely sturdy $\Tcontr$-presheaves, then 
\[
\begin{tikzcd}
\nondeg A \rar \dar \drar["\lrcorner" very near start, phantom] & \lilh^* A  \dar\\
\nondeg B \rar & \lilh^* B
\end{tikzcd}
\]
is a pullback by \cref{lem cartesian unit}.
\end{remark}

\begin{corollary}\label{cor JX identification}
  If $X\in \pre(\Tsimpl)$ is completely sturdy then $\lilk_!\lilk^* X \simeq \lilj_!
  \nondeg X$.
\end{corollary}

\begin{proof}
For $X$ completely sturdy, by definition 
  $\nondeg X$ is the unique 
  solution to the equation $\lilh_! \nondeg X = \lilk^* X$, so it follows that 
$\lilj_! \nondeg X = \lilk_! \lilh_! \nondeg X = \lilk_! \lilk^* X$.
\end{proof}

If $X\colon \simpl^\op \to \spaces$ is a stiff and complete simplicial space, then we also have $\lilk_! \lilk^* X \simeq \lilj_! \nondeg X$ by applying the corollary to $\lilp^* X$.

\begin{remark}
  The assignment $X \mapsto \vec{X}$ is functorial only in conservative maps between
  completely sturdy objects in $\pre(\Tsimpl)$ (since it makes use of \cref{thm
  restricted slices}), hence so too is the assignment $X \mapsto \lilj_! \vec{X}$
  (which earlier we wrote as $J(X)$). We note, however, that $\lilk_!\lilk^*$ is defined
  on all $\Tsimpl$-presheaves and is functorial in all maps.
\end{remark}

\section{Ikeo maps and left Kan extension along \texorpdfstring{$\lilk$}{k}}

\label{sec:ikeo}

In this section we establish a few key properties of $\lilk_! \colon \pre(\Tcontr) \to \pre(\Tsimpl)$.
By \cref{lem k shriek at T}, we have $(\lilk_! A)_{\Tob} = A_{\Tob}$.
By \cref{beck-chevalley} we have \[ (\lilk_!A)_r  =  \colim_{\alpha: [r] \ract [n]} A_n\]  
where the colimit is over the opposite of the category $\simpl^- \times_{\simpl^\actrm}
\simpl^\actrm_{[r]/}$
whose objects are active maps $\alpha: [r] \ract [n]$ and whose morphisms
are triangles 
\[
\begin{tikzcd}[column sep={2em,between origins}, row sep={7ex,between origins}]
 & {[r]} \ar[ld, -act] \ar[rd, -act] &  \\
{[n]} \ar[rr, twoheadrightarrow] & & {[m]}
\end{tikzcd}
\]
(codegeneracy maps mediating between active maps). 
An object in this category is an active map $\alpha\colon 
[r] \ract [n]$, which is uniquely the amalgamated sum of $r$ active maps $\alpha_i \colon [1] 
\ract [n_i]$, and in particular $[n] = [n_1] \vee \cdots \vee [n_r]$. 
The assignment $\alpha \mapsto (\alpha_1, \dots, \alpha_r)$ induces an isomorphism of categories
\begin{equation}
\label{eq Dop decomp}
  \simpl^- \times_{\simpl^\actrm} \simpl^\actrm_{[r]/} 
  \ \cong \ 
  \prod_{u=1}^r 
  \big( \simpl^- \times_{\simpl^\actrm} \simpl^\actrm_{[1]/} \big) \cong \prod_{u=1}^r \simpl^-,
\end{equation}
since there is a unique active map $[1] \ract [n_i]$.
Note that $(\lilk_!A)_0 = A_0$.

\begin{lemma}\label{lem unit ikeo}
If $A \in \pre(\Tcontr)$ is completely sturdy, then the unit $A \to \lilk^* \lilk_! A$ is ikeo.
\end{lemma}
\begin{proof}
At $\Tob$, the map is an equivalence by \cref{lem k shriek at T}.
To show this map ikeo, we must thus show that the square 
\[ \begin{tikzcd}
A_r \rar \dar & A_1^{\times r}  \dar \\
(\lilk^* \lilk_! A)_r \rar & (\lilk^* \lilk_! A)_1^{\times r}
\end{tikzcd} \]
is a pullback for all $r\geq 0$. Since $A$ is completely sturdy, we have $A\simeq \lilh_! \nondeg A$, so this square is equivalent to 
\[ \begin{tikzcd}
(\lilh_! \nondeg A)_r \rar \dar & (\lilh_! \nondeg A)_1^{\times r} \dar 
    \\
(\lilk^* \lilj_! \nondeg A)_r \rar & (\lilk^* \lilj_! \nondeg A)_1^{\times r} .
\end{tikzcd} \]
Using the formulas for
$\lilh_!$ and $\lilj_!$ we now expand this 
square to
\[
\begin{tikzcd}[ds]
\sum_{[r] \twoheadrightarrow [k]} \nondeg A_k \rar \dar& 
\prod_{i=1}^r \sum_{[1] \twoheadrightarrow [k_i]} \nondeg A_{k_i} \rar{\cong} \dar&
\sum_{(\ob \simpl^-_{[1]/})^r} \nondeg A_{k_1} \times \cdots \times \nondeg A_{k_r} \dar
\\
\sum_{[r] \ract [m]} \nondeg A_m \rar & 
\prod_{i=1}^r \sum_{[1] \ract [m_i]} \nondeg A_{m_i} \rar{\cong} &
\sum_{(\ob \simpl^\actrm_{[1]/})^r} \nondeg A_{m_1} \times \cdots \times \nondeg A_{m_r} 
\end{tikzcd}
\]
where in the top row we are summing over surjections out of $[r]$, and in the bottom row we are summing over active maps out of $[r]$.
But an active map (resp.\ surjection) out of $[r]$ is the same thing as $r$ active maps (resp.\ surjections) out of $[1]$.
We thus can see the preceding square as a sum over all $[r] \ract [m]$ of the following squares: the left in the case when the active map is surjective, and the right in the case when it is not.
\[ \begin{tikzcd}
\nondeg A_m \rar \dar&
\nondeg A_{m_1} \times \cdots \times \nondeg A_{m_r} \dar
    &
\varnothing \rar \dar & \varnothing \dar 
\\
\nondeg A_m \rar &
\nondeg A_{m_1} \times \cdots \times \nondeg A_{m_r} 
    &
\nondeg A_m \rar &
\nondeg A_{m_1} \times \cdots \times \nondeg A_{m_r} 
\end{tikzcd} \]
Both of these squares are pullbacks, so we may conclude that the unit is ikeo.
\end{proof}

\begin{lemma}\label{lem 2-Segal adjunct ikeo}
Let $X$ be 2-Segal.
If $A \to \lilk^* X$ is ikeo, then so is its adjunct $\lilk_! A \to X$.
\end{lemma}
\begin{proof}
Given an active map $[r] \ract [m]$, each of the following squares are pullbacks:
\[ \begin{tikzcd}[ds]
A_m \dar \rar \drar[phantom, "\lrcorner" very near start] & X_m \dar \rar \drar[phantom, "\lrcorner" very near start] & X_r \dar \\
\prod_{i=1}^r A_{m_i} \rar & \prod_{i=1}^r X_{m_i} \rar & \prod_{i=1}^r X_1 ,
\end{tikzcd} \]
the left since our map is ikeo, the right since $X$ is 2-Segal by \cite[Proposition 6.9]{GKT1}.
We thus have an equivalence $A_m \simeq X_r \times_{X_1^{\times r}} (A_{m_1} \times \cdots \times A_{m_r})$,  natural in $\simpl^- \times_{\simpl^\actrm} \simpl^\actrm_{[r] /} \cong (\simpl^-)^{\times r}$.
Taking colimits, we have
\[
    (\lilk_! A)_r \simeq \colim A_m \simeq \colim \left( X_r \times_{X_1^{\times r}} \prod_{i=1}^r A_{m_i} \right) \simeq X_r \times_{X_1^{\times r}} \colim \left( \prod_{i=1}^r A_{m_i} \right).
\]
We actually have  $\colim \left( \prod_{i=1}^r A_{m_i} \right) \simeq (\lilk_! A)_1^{\times r}$.
This is because
\[
\colim_{((\simpl^-)^\op)^{\times r}} \prod_{i=1}^r A_{m_i} \simeq \colim_{(\simpl^-)^\op} \colim_{(\simpl^-)^\op} \cdots \colim_{(\simpl^-)^\op} \prod_{i=1}^r A_{m_i} \simeq \prod_{i=1}^r \colim_{(\simpl^-)^\op} A_{m_i} \simeq \prod_{u=1}^r \lilk_!(A)_1 ,
\]
the middle equivalence since colimits in $\basecat$ commute with products by a fixed object.
We conclude that $(\lilk_! A)_r \simeq X_r \times_{X_1^{\times r}} (\lilk_!A)_1^{\times r}$, so $\lilk_! A \to X$ is ikeo.
\end{proof}

\begin{proposition}\label{prop ikeo iff adjunct ikeo}
Suppose $A \in \pre(\Tcontr)$ is completely sturdy and $X\in\pre(\Tsimpl)$ is 2-Segal.
A map $A \to \lilk^* X$ is ikeo if and only if its adjunct $\lilk_! A \to X$ is so.
\end{proposition}
\begin{proof}
The forward direction follows from \cref{lem 2-Segal adjunct ikeo}.
For the reverse direction, we use that both maps in the composite 
\[
    A \to \lilk^* \lilk_! A \to \lilk^* X
\]
are ikeo maps, by \cref{lem unit ikeo} and the fact that $\lilk^*$ preserves ikeo maps.
\end{proof}

\begin{lemma}\label{new lem ikeo and k shriek}
Let $f\colon B\to A$ be a map in $\pre(\Tcontr)$ with $A$ and $B$ completely sturdy.
Then $f$ is ikeo if and only if $\lilk_!f$ is ikeo in $\pre(\Tsimpl)$.
\end{lemma}
\begin{proof}
Notice that $\lilk_! A = \lilj_! \nondeg A$ is 2-Segal.
By \cref{prop ikeo iff adjunct ikeo}, $\lilk_! f$ is ikeo if and only if 
\[ \begin{tikzcd}
B \rar{f} & A \rar{\eta_A} & \lilk^* \lilk_! A
\end{tikzcd} \]
is ikeo.
As $\eta_A$ is ikeo by \cref{lem unit ikeo}, we have $f$ is ikeo if and only if this composite is ikeo. (Composition + cancellation of ikeo maps, which both follow from the pasting law for pullbacks.)
\end{proof}

\begin{lemma}\label{lem ikeo counit}
If $X \in \pre(\Tsimpl)$ is 2-Segal, then the counit $\varepsilon_X \colon \lilk_! \lilk^*X \to X$
is ikeo.
\end{lemma}
\begin{proof}
Apply \cref{lem 2-Segal adjunct ikeo} to the identity map $\lilk^* X \to \lilk^*X$.
\end{proof}

The following will be used in the next section only to see that $J(M) \to M$ is actually an object of the category $\dgr{M}$ we wish to show it is terminal in.

\begin{lemma}\label{counit conservative}
  If $M\in \pre(\Tsimpl)$ is 2-Segal 
  with indecomposable identities, then $\varepsilon_M \colon \lilk_!\lilk^*M \to M$ is
  conservative.
\end{lemma}

\begin{proof}
First, we observe that the inclusion $I \colon \simpl^-_{>0} \hookrightarrow \simpl^-$ is initial.
We must show that $I\downarrow [n]$ is contractible for all $n\geq 0$. 
This is clear for $n > 0$, since the comma category contains the terminal object $\id_{[n]}$.
For $n=0$, we have $I \downarrow [0] \cong \simpl^-_{>0}$.
Under Joyal duality, the opposite of this category is isomorphic to $\simpl^+$, hence is contractible (cf.\ the $n=0$ case of the proof of \cite[\href{https://kerodon.net/tag/04RE}{Tag 04RE}]{kerodon}).

Since $M$ has indecomposable identities, we know
\[ \begin{tikzcd}
M_0 \rar \drar[phantom, "\lrcorner" very near start] \dar & M_n \dar \\
M_0 \rar & M_1
\end{tikzcd} \]
is a pullback for all $n > 0$.
Finality says that $(\lilk_! \lilk^* M)_1 = \colim_{(\simpl^-)^\op} M_n \simeq \colim_{(\simpl^-_{>0})^\op} M_n$, so we have
\[
  M_0 \times_{M_1} (\lilk_! \lilk^* M)_1 \simeq 
  M_0 \times_{M_1} \colim_{(\simpl^-_{>0})^\op} M_n 
\simeq
\colim_{(\simpl^-_{>0})^\op} (M_0 \times_{M_1} M_n) \simeq \colim_{(\simpl^-_{>0})^\op} M_0.
\]
Since $(\simpl^-_{>0})^\op$ is contractible, this last term is equivalent to $M_0$.
We conclude that 
\[ \begin{tikzcd}
M_0 \rar \dar & (\lilk_! \lilk^* M)_1 \dar 
\ar[drr,"\simeq" description, phantom]
& & (\lilk_! \lilk^* M)_0 \rar \dar & (\lilk_! \lilk^* M)_1 \dar\\
M_0 \rar & M_1 
& & 
M_0 \rar & M_1 
\end{tikzcd} \]
is a pullback.
Combining with \cref{lem ikeo counit}, we see that the counit is conservative.
\end{proof}

\section{The main theorem}\label{sec main thm proof}

\begin{theorem}\label{main theorem}
Suppose $M \in \pre(\Tsimpl)$ is 2-Segal with indecomposable identities, standard, and completely sturdy.
Then $J(M) \to M$ is a terminal object in $\dgr{M}$,  the category of conservatively graded 2-Segal presheaves over $M$.
\end{theorem}
\begin{proof}
This is an object of $\dgr{M}$ by \cref{counit conservative}.
It is terminal by \cref{main lemma} below.
\end{proof}

Notice that under these hypotheses, $J(M)$ possesses the special zeta function (cf.\ \ref{zeta of free}), since we have $(\lilk_! \lilk^* M)_\Tob = M_0 + \nondeg M_\Tob = M_0 + \nondeg M_1$ (by \cref{cor JX identification} and \cref{lke of hj}).
The classical case is when $M$ is the $\Tsimpl$-presheaf given as
$B\N$ with standard zeta function, yielding \cref{intro theorem A}.

\begin{lemma}\label{lem slice pullback}
Suppose 
$L : \mathcal{C} \rightleftarrows \mathcal{D} : R$ is an adjunction, and consider for each $d\in D$
the induced adjunction $L_d : \mathcal{C}_{/Rd} \rightleftarrows \mathcal{D}_{/d} : R$.
Then the square below left is a pullback.
\[ \begin{tikzcd}
\mathcal{C}_{/Rd} \rar{L_d} \dar \drpullback & \mathcal{D}_{/d} \dar & & \mathcal{C}_{/Rd} \rar{L_d} \dar \drpullback & \mathcal{D}_{/d} \dar \\
\mathcal{C} \rar[swap]{L} & \mathcal{D} & & \mathcal{C}_{/Rd'} \rar[swap]{L_{d'}}  & \mathcal{D}_{/d'}
\end{tikzcd} 
\]
It follows that for each $d\to d'$ in $\mathcal{D}$, the square above right is a pullback.
\end{lemma}
\begin{proof}
The first statement is the dual of \cite[Lemma 2.1.1]{SchlankYanovski:ICEHA}.
The second follows from the pasting law for pullbacks.
\end{proof}

\begin{lemma}\label{main lemma}
  Let $M \in \pre(\Tsimpl)$ be standard and completely sturdy, and $\lilk_!\lilk^* M = J(M) \to M$ the counit. If $X \in \pre(\Tsimpl)$ is 2-Segal and
  $X\to M$ is conservative, then the category
  \[
    \pre(\Tsimpl)^\ikeo_{/X} \underset{\pre(\Tsimpl)_{/M}}\times \pre(\Tsimpl)^\culf_{/J(M)}
  \]
  is a contractible groupoid.
\end{lemma}
\begin{proof}
In the following cube whose front and back squares are pullbacks, we also have that the bottom square is a pullback by \cref{lem slice pullback}.
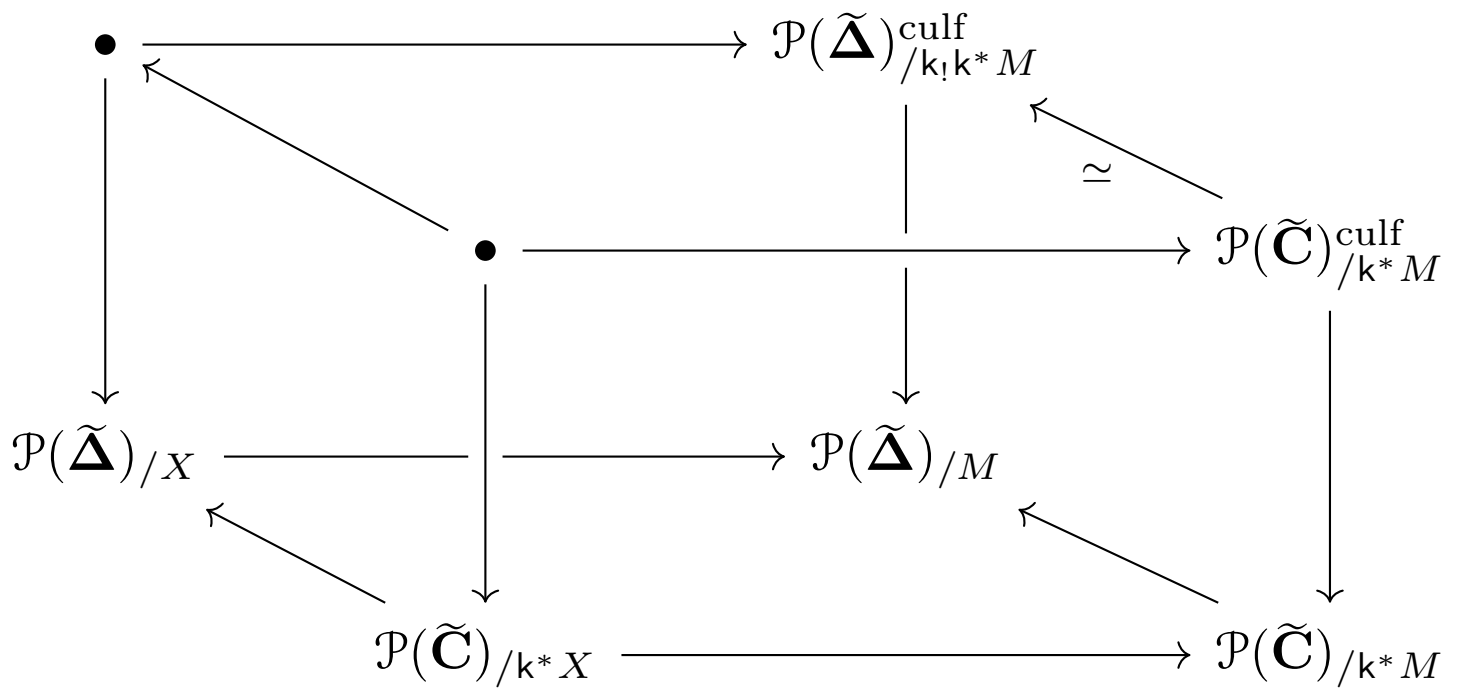
\[ \begin{tikzcd}
\bullet\ar[rr] \ar[dd] & & \pre(\Tsimpl)^\culf_{/\lilk_!\lilk^*M} \ar[dd] \\
& \bullet \ar[rr, crossing over] \ar[ul] & &  \pre(\Tcontr)^\culf_{/\lilk^*M}  \ar[dd]  \ar[ul,"\simeq"] \\
\pre(\Tsimpl)_{/X} \ar[rr] & & \pre(\Tsimpl)_{/M} \\
& \pre(\Tcontr)_{/\lilk^* X} \ar[rr] \ar[ul] \ar[from=uu, crossing over] & & \pre(\Tcontr)_{/\lilk^* M} \ar[ul]
\end{tikzcd} \]
Therefore the top square is a pullback.
Since $M$ is completely sturdy, the restriction $\lilk^*M$ is in the image of $\lilh_!$, so \cref{thm four equivs} implies that the upper right diagonal edge is an equivalence. 
Hence the upper left diagonal edge is an equivalence as well.

We now restrict to full subcategories as displayed below.
\begin{equation}\label{diagram restrict to ikeo}
\begin{tikzcd}
\displaystyle \pre(\Tcontr)^\ikeo_{/\lilk^* X}  \underset{\pre(\Tcontr)_{/\lilk^* M}}\times \pre(\Tcontr)^\culf_{/\lilk^*M}
\rar[hook]
\dar 
&
\displaystyle \pre(\Tcontr)_{/\lilk^* X}  \underset{\pre(\Tcontr)_{/\lilk^* M}}\times \pre(\Tcontr)^\culf_{/\lilk^*M}
\dar{\simeq}
 \\
\displaystyle \pre(\Tsimpl)^\ikeo_{/X}  \underset{\pre(\Tsimpl)_{/M}}\times \pre(\Tsimpl)^\culf_{/\lilk_!\lilk^*M} 
\rar[hook]
&
\displaystyle \pre(\Tsimpl)_{/X}  \underset{\pre(\Tsimpl)_{/M}}\times \pre(\Tsimpl)^\culf_{/\lilk_!\lilk^*M}
\end{tikzcd} 
\end{equation} 
This requires justification, as it is not clear a priori that
$\pre(\Tcontr)^\ikeo_{/\lilk^* X} \to \pre(\Tcontr)_{/\lilk^* X} \to
\pre(\Tsimpl)_{/X}$ takes values in $\pre(\Tsimpl)^\ikeo_{/X}$. But the left map
of \eqref{diagram restrict to ikeo} \emph{does} land in the indicated full
subcategory: given $A \to \lilk^* X$ such that the composite $A \to \lilk^* X
\to \lilk^* M$ is culf, since $M$ is completely sturdy we deduce that $A$ is as well.
Since $X$ is 2-Segal 
and $A$ is completely sturdy, $A \to \lilk^*X$ is ikeo if and only if $\lilk_!A \to X$ is
ikeo (\cref{prop ikeo iff adjunct ikeo}).

Now suppose we have an object
$(X\xleftarrow{\ikeo} U \xrightarrow{\culf} \lilk_!\lilk^* M)$ in the bottom left of \eqref{diagram restrict to ikeo}.
By the right-hand equivalence, we automatically have \[ (U \to \lilk_! \lilk^* M) \simeq \lilk_!( A  \to \lilk^* M)\]
for some $A\to \lilk^* M$, which then must be the composite $A \to \lilk^* X \to \lilk^*M$. 
But then \eqref{diagram restrict to ikeo} is a pullback by \cref{prop ikeo iff adjunct ikeo}, so the left vertical edge is an equivalence. 

To prove \cref{main theorem}, it remains to show that
\begin{equation}\label{eq contr reduction}
    \pre(\Tcontr)^\ikeo_{/\lilk^* X}  \underset{\pre(\Tcontr)_{/\lilk^* M}}\times \pre(\Tcontr)^\culf_{/\lilk^*M}
\end{equation}
is contractible. This is the full subcategory of $\pre(\Tcontr)^\ikeo_{/\lilk^*
X}$ on those ikeo maps $A \to \lilk^* X$ such that the composite $A \to \lilk^*X
\to \lilk^*M$ is culf. Since $M$ is standard and this composite is culf, also $A$ is
standard. This tells us that \eqref{eq contr reduction} is a full
subcategory of $\pre(\Tcontr)^{\ikeo,\std}_{/\lilk^* X}$ consisting of ikeo maps
with standard domain. In fact, \eqref{eq contr reduction} is the category
$\pre(\Tcontr)^{\ikeo,\std,\cons}_{/\lilk^* X}$, because if $A$ is standard and $A
\to \lilk^* X$ is conservative, then since $X \to M$ is conservative, the
composite $A \to \lilk^*M$ is culf (\cref{rmk culf vs conservative}); conversely if $A \to \lilk^*M$ is
conservative, then cancellation shows that $A\to \lilk^*X$ is conservative. 
But
\cref{prop std ikeo contractible} says that $\pre(\Tcontr)^{\ikeo,\std,\cons}_{/\lilk^* X}$
is a contractible groupoid, so the pullback category \eqref{eq contr reduction} is as well, as we wanted to show.
\end{proof}

\begin{blanko}[Unraveling \texorpdfstring{$W$}{W}]\label{unravel}
  Having given a rather abstract proof of uniqueness of 
  the span $X\leftarrow W \to 
  \lilk_!\lilk^*M$,
  we owe a description of the span in more explicit terms.
With reference to the proof of \cref{main theorem} we consider first
the standard $\Tcontr$-presheaf $A$ with its ikeo and conservative map $A \to \lilk^* X$.
The left leg of the span is the adjunct map $\lilk_! A \to X$.
Since $A$ is standard and the map $A \to \lilk^* X$ is ikeo, we have equivalences
\begin{equation}\label{A1=XT}
  \begin{tikzcd}[row sep=0pt, baseline=5.2pt]
A_1 & A_\Tob \lar["\sim"'] \rar["\sim"] & \lilk^*X_\Tob = X_\Tob \\
& A_0 \rar["\sim"] & \lilk^*X_0 = X_0.
\end{tikzcd}
\end{equation}
Since $\lilk^*M$ is completely sturdy and 
the two maps $A\to\lilk^* X\to \lilk^*M$ are conservative,
also the presheaves $A$ and $\lilk^*X$ are completely sturdy, so $A \simeq \lilh_! \vec{A}$ and
$\lilk^* X \simeq \lilh_! \vec{X}$.

Since $A \to \lilk^* X$ is conservative, we have by \cref{rmk functoriality of nondeg}
the pullback square on the left
\[ \begin{tikzcd}
\vec{A}_r \rar \dar \drar[phantom, "\lrcorner" very near start]& A_r \rar \dar 
\drar[phantom, "\lrcorner" very near start] & A_1^{r} \dar \\
\vec{X}_r \rar & X_r \rar & X_1^{r} ,
\end{tikzcd} \]
and since $A \to \lilk^* X$ is ikeo, also the right square is
a pullback.
From the composite pullback together with \eqref{A1=XT}
we see that
\[
  \vec{A}_r = \vec{X}_r \times_{X_1^r} X_\Tob^r.
\]
But now $W = \lilk_! A = \lilk_! \lilh_! \vec{A} = \lilj_! \vec{A}$, and we arrive at
the formula
\[
  W_n = \sum_{\alpha\colon [n]\ract[r]} \big( \nondeg X_r \times_{X_1^r} X_\Tob^r \big)
\]
by \cref{lke of hj}.
Since $W \to X$ is ikeo, we have $W_0 \simeq X_0$ and $W_\Tob \simeq X_\Tob$.

Note that in the particular case where $X$ is standard (so that $X_\Tob=X_1$), we thus
have $W = \lilk_! \lilk^* X = J(X)$ and the span is simply
$$
X \leftarrow J(X) \rightarrow J(M)  .
$$
Here $X\leftarrow J(X)$ is ikeo by \cref{lem ikeo counit} and 
$J(X)\to J(M)$ is culf by \cref{thm restricted slices}.
\end{blanko}

\section{Finiteness conditions}\label{sec finiteness}

All the results so far
belong 
to the objective level, where it is not necessary to impose any 
finiteness conditions,
and in Sections~\ref{sec object prop}--\ref{sec main thm proof}
we have even worked over a rather general category $\basecat$. To come back to
classical algebraic combinatorics in vector spaces, we now return to the case
$\basecat=\spaces$, and impose the necessary finiteness conditions to ensure that we
can take homotopy cardinality and obtain results over vector spaces.
We'll see that the finiteness conditions required
for the input data to admit a cardinality (and thereby admit incidence coalgebras 
in the category of $\Q$-vector spaces) play well together with the constructions and
the universal property: if we impose these conditions on the input ($M$ and $X$) then 
also the output ($J(M)$ and the span $X\leftarrow W \to J(M)$) will satisfy these conditions.

The finiteness conditions required to take cardinality are
generally that all spans should be of finite type \cite{Galvez-Kock-Tonks:1602.05082},
meaning that the backward map is a 
finite map. For the decomposition space $X$ itself, this condition must be imposed on 
the spans that define the coalgebra structure, so $X_1 \stackrel{d_1}\leftarrow X_2$ and
$X_1 \stackrel{s_0}\leftarrow X_0$ should be required to be finite maps,
so that the sum in the formula
for comultiplication
$$
\Delta(\xora) = \sum_{\substack{\sigma\in X_2 \\ d_1(\sigma)=\xora}} d_2(\sigma) \otimes 
d_0(\sigma)$$
is a finite sum. 
Equivalently, all active 
maps of $X$ should be finite, which is to say it is a locally finite decomposition 
space~\cite{GKT2}. (Note that the condition on $s_0$ is already implied by completeness.)

For the zeta function to have a cardinality we need $X_1 \leftarrow X_\Tob$ to be a 
finite map. Since $\sigr^*$ is active, the natural generalization of local finiteness 
to $\Tsimpl$-presheaves is therefore the following.

\begin{definition}
Let $X$ be a $\Tsimpl$-presheaf whose underlying simplicial space $\lili^*X$ is a complete decomposition space.
Then $X$ is called
\begin{itemize}
\item \emph{locally finite} when all its active maps are finite (equivalently, its underlying simplicial space is locally finite, and the map $X_1 \leftarrow X_\Tob$ is finite);
\item \emph{of locally finite length} when $\lili^*X$ is so;
\item \emph{M\"obius} if it is sturdy, locally finite, and of locally finite length.
\end{itemize}
\end{definition}

\begin{lemma}\label{J finite}
  For any sturdy $\Tsimpl$-presheaf $X$ with finite 
  $X_1 \leftarrow X_\Tob$ (for example $X$ standard), the $\Tsimpl$-presheaf $J(X) \coloneq
  \lilk_! \lilk^* X = \lilj_! \nondeg X$ is M\"obius. (No finiteness conditions
  are required on the underlying simplicial space of $X$.)
\end{lemma}
\begin{proof}
  That the underlying simplicial space of $J(X) = \lilj_! \nondeg X$ is M\"obius is 
  the content of 
  \cite[Corollary 2.3.3]{Hackney-Kock:2210.11192}. It remains to check that the active
  map $J(X)_1 \leftarrow J(X)_\Tob$ is finite. But this map is the composite
  $$
  \sum_r \nondeg X_r \leftarrow X_0+\nondeg X_1 = X_1 \leftarrow X_\Tob
  $$
  where the rightmost map is finite by assumption and the leftmost map is finite 
  because it is a sum inclusion.
\end{proof}

\begin{proposition}\label{W finite}
Suppose $M$ and $X$ are M\"obius $\Tsimpl$-presheaves and $M$ is standard.
If $X$ is conservatively graded over $M$, then the universal span $X\leftarrow W \to J(M)$ is of finite type.
\end{proposition}
\begin{proof}
By Corollary 5.11 and Proposition 6.6 of \cite{GKT2}, $M$ has indecomposable identities; it thus satisfies the hypotheses of \cref{main theorem}.
The map
  $X_1 \leftarrow W_1$ is the composite
  $$
  X_1 \leftarrow \sum_r \nondeg X_r \leftarrow \sum_r \big( \nondeg X_r
  \times_{X_1^r} X_\Tob^r \big) .
  $$
  Since $X$ is assumed M\"obius, the first map is finite 
  by \cref{lambda finite} below, and the second is finite
  because 
  $X_1 \leftarrow X_\Tob$ is finite (by finiteness assumption on $X_\Tob$).
\end{proof}

\begin{lemma}\label{lambda finite}
  Let $X \in \pre(\simpl)$ be a complete decomposition space. Then $X_1 \leftarrow
  J(X)_1$ is a finite map if and only if $X$ is M\"obius.
\end{lemma}
\begin{proof}
  Recall that M\"obius means locally finite and of locally finite length.
  The map in question is the ``long-edge map'' $X_1 \leftarrow \sum_r \nondeg X_r$.
  To ask that the sum is finite for each fixed $1$-simplex $\xora\in 
  X_1$ is to say that $X$ is of locally finite length. To ask furthermore that for 
  each $r$ the map $X_1 \leftarrow \nondeg X_r$ is finite is equivalent to asking that
  $X$ is locally finite.
\end{proof}

\begin{remark}
  It is common in algebraic combinatorics to impose the condition
  on graded (co)algebras that the graded pieces are finite dimensional. (This 
  condition is
  not required for the universal property of $\QSym$, though.) In our situation we
  consider more general notions of grading, but it is still a natural condition to
  impose that the graded pieces be finite dimensional, which at the objective level
  means precisely that the conservative map expressing the grading is a finite map. The grading
  of the terminal object $J(M)$ of the main theorem is given by the map
  $$
  M \leftarrow J(M) .
  $$
  If we want this map to be finite, then \cref{lambda finite} tells us we should
  demand that $M$ be M\"obius. And then $X$ is automatically M\"obius too (by
  \cref{lem:locally finite length}) provided it is locally finite, and 
  ``locally finite'' 
  in turn is automatic if we also ask the grading of $X$ to be a finite map. In
  conclusion, a natural setting for the finite version of the main theorem is to
  stipulate that a grading should be a (conservative) finite map with values in a
  M\"obius decomposition space (with standard zeta function).
  
  With this assumption on gradings, there are {\em no} further conditions to impose on $X$.
\end{remark}

\begin{blanko}[The Aguiar--Bergeron--Sottile formula]
  We finally consider the classical case, where $M=B\N$ and $J(M)=Q$ and show 
  how the Aguiar--Bergeron--Sottile formula drops out of our general construction.
  Instantiating here the explicit formula from  \ref{unravel}
  for a general $(X,F)$ (M\"obius, with zeta function with finite fibres, that is, 
  finite values),
  we find that the universal span in the classical case is
  $$
  X_1 \longleftarrow  \sum_r \big( \nondeg X_r \times_{X_1^r} F^r \big) 
  \longrightarrow 
  \sum_r \N_+^r   .
  $$
  The left leg is ``long edge''; the right leg is ``degree vector of the principal 
  edges''.
  The fact that there are only nondegenerate simplices in the left factor of 
  the
  summands of
  the apex (that is, $\nondeg X_r$ rather than $X_r$) can also be expressed in the right factor: we have
  $$
  \nondeg X_r \times_{X_1^r} F^r \simeq X_r \times_{X_1^r} (F\backslash X_0)^r
  $$
  so that the span takes the form
  $$
  X_1 \longleftarrow  \sum_r \big( X_r \times_{X_1^r} (F\backslash X_0)^r \big) 
  \longrightarrow  
  \sum_r \N_+^r   .
  $$
  We can now read off how the span acts on some element $\name x : 1 \to X_1$ in the 
  incidence coalgebra of $X$:
  we sum over all $r$ and then sum over all $r$-simplices with long edge $x$.
  This is to sum over all $r$-ary comultiplications. On each such $r$-simplex, we 
  extract the $r$ factors (the $r$ principal edges of the $r$-simplex) and evaluate $F\backslash X_0$ on 
  each factor. Finally the right leg of the span takes each such principal edge $x_i$ 
  to its degree, and altogether returns the degree vector as an element of $\N_+^r$;
  each corresponds to monomial quasisymmetric function. In conclusion, if $(C,\phi)$ is
  the incidence coalgebra of $(X,F)$, meaning that $\phi$ is the cardinality of $F$ 
  (and $\phi-\epsilon$ is the cardinality of $F\backslash X_0$), then we see that
  the linear map $C\to \QSym$ implemented by the 
  span is precisely the Aguiar--Bergeron--Sottile universal map
  $$
  x \mapsto \sum_{r\geq 0} \sum_{x_{(1)}\otimes\dots\otimes x_{(r)}\in \Delta^{(r)}(x)}
  (\phi-\epsilon)^{\otimes r}(x_{(1)}\otimes\dots\otimes x_{(r)}) \cdot 
  M_{\deg(x_{(1)})\cdots \deg(x_{(r)})} 
  $$
  discussed in \ref{ABS}.
  (We do not see from this formula that it is a coalgebra map, but that we know 
  already from the fact that the span is ikeo-culf.)
\end{blanko}

\appendix

\section{Subdividing a morphism is \texorpdfstring{$\infty$}{∞}-categorical: \texorpdfstring{$\Tsimpl$}{Δ̃} as a pushout}\label{app subdividing}

In \cref{sec shape categories} we described $\Tsimpl$ and $\Tcontr$ by pushouts of 1-categories $\simpl \amalg_{[1]} [2]$ and $\contr \amalg_{[1]} [2]$.
These pushouts may instead be taken in $\catinf$, a special case of the following proposition.
This may be regarded as a variation of the situation in \cite[\S4.3]{HORR:I2CPT}.
We use techniques for analyzing pushouts in $\cat$ and $\catinf$ from \cite{HORR-Dwyer,HaineRamziSteinebrunner:FFFPO}. 

\begin{proposition}\label{prop inf pushout}
Suppose $C$ is a small category, $f \colon x \to y$ is a morphism in $C$, and $n\geq 1$.
Let
\[ \begin{tikzcd}
{[1]} \dar[-act] \rar{\name{f}} \ar[dr, phantom, "\ulcorner" very near end] & C \dar \\
{[n]} \rar & D
\end{tikzcd} \]
be the pushout in $\cat$.
Then this square is also a pushout in $\catinf$.
\end{proposition}
\begin{proof}
In \cref{lem discretely flat} below we show that $[1] \ract [n]$ is a discretely flat inclusion, so the result follows from Theorem 1.5 in \S4 of \cite{HORR-Dwyer}.
\end{proof}

\begin{lemma}\label{lem pushout calc}
Assume the same setup as \cref{prop inf pushout}, and for $0 \leq i \leq j \leq n$ write $\gamma_{ij} \colon i \to j$ for the unique morphism in $[n]$.
Then the hom sets of $D$ are described as follows, where $p,q\in C$:
\begin{itemize}
\item $D(p,q) = C(p,q)$
\item $D(k,q) \cong C(y,q)$ via $g\circ \gamma_{kn} \mapsfrom g$ for $0 < k < n$
\item $D(p,k) \cong C(p,x)$ via $\gamma_{0k} \circ g \mapsfrom g$ for $0 < k < n$
\item $D(i,j) \cong C(y,x) \amalg \{ \gamma_{ij} \}$ via $\gamma_{0j} \circ g \circ \gamma_{in} \mapsfrom g$ for $0 < i \leq j < n$ and $g\colon y \to x$, and 
\item $D(i,j) \cong C(y,x)$ via $\gamma_{0j} \circ g \circ \gamma_{in} \mapsfrom g$ for $0 < j < i < n$. \qed
\end{itemize}
\end{lemma}

\begin{lemma}\label{lem discretely flat}
The active map $[1] \ract [n]$ is a \emph{discretely flat inclusion} in the sense of Definition 3.4 of \cite{HORR-Dwyer}.
\end{lemma}
\begin{proof}
Let $\scat$ denote the category of simplicially-enriched categories equipped with the 
Bergner model structure, and $\disc \colon \cat \to \scat$ the functor that takes a 1-category to the simplicial category with discrete hom simplicial sets.
A discretely flat inclusion is first of all an inclusion, meaning that it is injective on objects and faithful, and which is discretely flat, in the sense that the pushout functor
\[\disc [n]  \amalg_{\disc [1] } (-) \colon \slicel{\disc [1]}{\scat} \to \slicel{\disc [n]}{\scat} \] 
preserves Dwyer--Kan equivalences.
The formulas in \cref{lem pushout calc} upgrade to isomorphisms of hom simplicial sets when taking the pushout $D$ of $C \leftarrow \disc [1] \rightarrow \disc [n]$.
Hence if $F \colon C \to C'$ is a Dwyer--Kan equivalence under $\disc [1]$, then $D \to D'$ is a Dwyer--Kan equivalence under $\disc [n]$.
Indeed, $D\to D'$ is homotopically fully-faithful by \cref{lem pushout calc} and 2-of-3, and the map of categories $\pi_0 D \to \pi_0 D'$ is essentially surjective since $\pi_0 C \to \pi_0 C'$ is so.
\end{proof}

\bibliographystyle{amsalpha}
\bibliography{universal}

\providecommand{\bysame}{\leavevmode\hbox to3em{\hrulefill}\thinspace}
\providecommand{\MR}{\relax\ifhmode\unskip\space\fi MR }
\providecommand{\MRhref}[2]{%
  \href{http://www.ams.org/mathscinet-getitem?mr=#1}{#2}
}
\providecommand{\href}[2]{#2}
\begin{thebibliography}{BMSvW00}

\bibitem[ABS06]{Aguiar-Bergeron-Sottile}
Marcelo Aguiar, Nantel Bergeron, and Frank Sottile, \emph{Combinatorial {H}opf
  algebras and generalized {D}ehn-{S}ommerville relations}, Compos. Math.
  \textbf{142} (2006), no.~1, 1--30. \MR{2196760}

\bibitem[AS21]{Alexandersson-Sulzgruber:1807.02460}
Per Alexandersson and Robin Sulzgruber, \emph{P-partitions and p-positivity},
  Int. Math. Res. Not. \textbf{2021} (2021), no.~14, 10848--10907.

\bibitem[Bar17]{Barwick:SMFEAK1}
Clark Barwick, \emph{Spectral {M}ackey functors and equivariant algebraic
  {$K$}-theory ({I})}, Adv. Math. \textbf{304} (2017), 646--727. \MR{3558219}

\bibitem[BB11]{Billera-Brenti:0710.3965}
Louis~J. Billera and Francesco Brenti, \emph{Quasisymmetric functions and
  {K}azhdan--{L}usztig polynomials}, Isr. J. Math. \textbf{184} (2011),
  317--348.

\bibitem[BHS25]{BarkanHaugsengSteinebrunner:EAP}
Shaul Barkan, Rune Haugseng, and Jan Steinebrunner, \emph{Envelopes for
  algebraic patterns}, Algebr. Geom. Topol. \textbf{25} (2025), no.~9,
  5319--5388. \MR{5005649}

\bibitem[BMSvW00]{Bergeron-Mykytiuk-Sottile-vanWilligenburg}
Nantel Bergeron, Stefan Mykytiuk, Frank Sottile, and Stephanie van
  Willigenburg, \emph{Noncommutative {P}ieri operators on posets}, J. Combin.
  Theory Ser. A \textbf{91} (2000), no.~1, 84--110.

\bibitem[BR08]{Baker-Richter:0605743}
Andrew Baker and Birgit Richter, \emph{Quasisymmetric functions from a
  topological point of view}, Math. Scand. \textbf{103} (2008), 208--242.

\bibitem[CH21]{Chu-Haugseng:1907.03977}
Hongyi Chu and Rune Haugseng, \emph{Homotopy-coherent algebra via {S}egal
  conditions}, Adv. Math. \textbf{385} (2021), Paper No. 107733, 95.
  \MR{4256131}

\bibitem[CK20]{Carlier-Kock:1807.11858}
Louis Carlier and Joachim Kock, \emph{Antipodes of monoidal decomposition
  spaces}, Commun. Contemp. Math. \textbf{22} (2020), no.~2, 1850081, 15.
  \MR{4077086}

\bibitem[DEFT20]{Diehl-EbrahimiFard-Tapia:1906.05823}
Joscha Diehl, Kurusch Ebrahimi-Fard, and Nikolas Tapia, \emph{Time-warping
  invariants of multidimensional time series}, Acta Appl. Math. \textbf{170}
  (2020), 265--290.

\bibitem[DK19]{DyckerhoffKapranov:HSS}
Tobias Dyckerhoff and Mikhail Kapranov, \emph{Higher {S}egal spaces}, Lecture
  Notes in Mathematics, vol. 2244, Springer, Cham, 2019. \MR{3970975}

\bibitem[Ehr96]{Ehrenborg}
Richard Ehrenborg, \emph{On posets and {H}opf algebras}, Adv. Math.
  \textbf{119} (1996), no.~1, 1--25. \MR{1383883}

\bibitem[GCKT]{GKT:QSym}
Imma G\'{a}lvez-Carrillo, Joachim Kock, and Andrew Tonks, \emph{Decomposition
  spaces of quasi-symmetric functions}, unpublished/in preparation.

\bibitem[GCKT18a]{GKT1}
\bysame, \emph{Decomposition spaces, incidence algebras and {M}\"{o}bius
  inversion {I}: {B}asic theory}, Adv. Math. \textbf{331} (2018), 952--1015.
  \MR{3804694}

\bibitem[GCKT18b]{GKT2}
\bysame, \emph{Decomposition spaces, incidence algebras and {M}\"{o}bius
  inversion {II}: {C}ompleteness, length filtration, and finiteness}, Adv.
  Math. \textbf{333} (2018), 1242--1292. \MR{3818099}

\bibitem[GCKT18c]{GKT3}
\bysame, \emph{Decomposition spaces, incidence algebras and {M}\"{o}bius
  inversion {III}: {T}he decomposition space of {M}\"{o}bius intervals}, Adv.
  Math. \textbf{334} (2018), 544--584. \MR{3828744}

\bibitem[GCKT18d]{Galvez-Kock-Tonks:1602.05082}
\bysame, \emph{Homotopy linear algebra}, Proc. Royal Soc. Edinburgh A
  \textbf{148} (2018), no.~2, 293--325. \MR{3777576}

\bibitem[GCKT24]{GKT:Crapo}
\bysame, \emph{Convex decomposition spaces and {C}rapo complementation
  formula}, 2024, \href{https://arxiv.org/abs/2409.03742}{arXiv:2409.03742}
  [math.CT].

\bibitem[GCKT26]{Galvez-Kock-Tonks:1612.09225}
\bysame, \emph{Decomposition spaces in combinatorics}, Higher {S}egal spaces
  and applications (J.~Bergner, J.~Kock, and M.~Sarazola, eds.), Contemp.
  Math., vol. 838, Amer. Math. Soc., Providence, RI, 2026, pp.~255--331.
  \MR{5061611}

\bibitem[Ges84]{Gessel:84}
Ira~M. Gessel, \emph{Multipartite {$P$}-partitions and inner products of skew
  {S}chur functions}, Combinatorics and algebra ({B}oulder, {C}olo., 1983),
  Contemp. Math., vol.~34, Amer. Math. Soc., Providence, RI, 1984,
  pp.~289--317. \MR{777705}

\bibitem[GHL22]{GagnaHarpazLanari:OEAMI2C}
Andrea Gagna, Yonatan Harpaz, and Edoardo Lanari, \emph{On the equivalence of
  all models for {$(\infty,2)$}-categories}, J. Lond. Math. Soc. (2)
  \textbf{106} (2022), no.~3, 1920--1982. \MR{4498545}

\bibitem[G{\"o}d24]{Godicke:IC2SS}
Jonte G{\"o}dicke, \emph{An $\infty$-category of 2-{S}egal spaces}, 2024,
  \href{https://arxiv.org/abs/2407.13357}{arXiv:2407.13357} [math.AT].

\bibitem[GR14]{Grinberg-Reiner:1409.8356}
Darij Grinberg and Victor Reiner, \emph{Hopf algebras in combinatorics}, 2014,
  \href{https://arxiv.org/abs/1409.8356}{arXiv:1409.8356} [math.CO].

\bibitem[Haz03]{Hazewinkel:2003}
Michiel Hazewinkel, \emph{Symmetric functions, noncommutative symmetric
  functions, and quasisymmetric functions}, Acta Appl. Math. \textbf{75}
  (2003), 55--83, Monodromy and diﬀerential equations (Moscow, 2001).

\bibitem[HHLN23]{HaugsengHebestreitLinskensNuiten:2VFFSCS}
Rune Haugseng, Fabian Hebestreit, Sil Linskens, and Joost Nuiten,
  \emph{Two-variable fibrations, factorisation systems and
  {$\infty$}-categories of spans}, Forum Math. Sigma \textbf{11} (2023), Paper
  No. e111, 70. \MR{4676218}

\bibitem[HK25]{Hackney-Kock:2210.11192}
Philip Hackney and Joachim Kock, \emph{Free decomposition spaces}, Collect.
  Math. \textbf{76} (2025), no.~3, 575--599. \MR{4950314}

\bibitem[HK26]{Hackney-Kock:2210.11191}
\bysame, \emph{Culf maps and edgewise subdivision}, Trans. Amer. Math. Soc.
  \textbf{379} (2026), no.~4, 2955--3007, with an appendix coauthored with Jan
  Steinebrunner. \MR{5039469}

\bibitem[Hof00]{Hoffman:9907173}
Michael~E. Hoffman, \emph{Quasi-shuffle products}, J. Algebr. Combin.
  \textbf{11} (2000), 49--68.

\bibitem[HORR23]{HORR:I2CPT}
Philip Hackney, Viktoriya Ozornova, Emily Riehl, and Martina Rovelli, \emph{An
  {$(\infty,2)$}-categorical pasting theorem}, Trans. Amer. Math. Soc.
  \textbf{376} (2023), no.~1, 555--597. \MR{4510118}

\bibitem[HORR24]{HORR-Dwyer}
\bysame, \emph{Pushouts of {D}wyer maps are {$(\infty,1)$}-categorical},
  Algebr. Geom. Topol. \textbf{24} (2024), no.~4, 2171--2183. \MR{4776277}

\bibitem[HRS25]{HaineRamziSteinebrunner:FFFPO}
Peter~J. Haine, Maxime Ramzi, and Jan Steinebrunner, \emph{Fully faithful
  functors and pushouts of $\infty$-categories}, 2025,
  \href{https://arxiv.org/abs/2503.03916}{arXiv:2503.03916} [math.CT].

\bibitem[Joy81]{Joyal:1981}
Andr{\'e} Joyal, \emph{Une th\'eorie combinatoire des s\'eries formelles}, Adv.
  Math. \textbf{42} (1981), no.~1, 1--82. \MR{MR633783 (84d:05025)}

\bibitem[Koc15]{Kock:1411.3098}
Joachim Kock, \emph{Perturbative renormalisation for not-quite-connected
  bialgebras}, Lett. Math. Phys. \textbf{105} (2015), no.~10, 1413--1425.
  \MR{3395225}

\bibitem[Koc20]{Kock:1809.00941}
\bysame, \emph{From {M}\"{o}bius inversion to renormalisation}, Commun. Number
  Theory Phys. \textbf{14} (2020), no.~1, 171--198. \MR{4060972}

\bibitem[KT97]{Krob-Thibon:NSym-IV}
Daniel Krob and Jean-Yves Thibon, \emph{Noncommutative symmetric functions.
  {IV.} {Q}uantum linear groups and {H}ecke algebras at q = 0}, J. Algebr.
  Combin. \textbf{6} (1997), 339--376.

\bibitem[Ler76]{Leroux:1976}
Pierre Leroux, \emph{Les cat{\'e}gories de {M}{\"o}bius}, Cahiers Topol.
  G{\'e}om. Diff. \textbf{16} (1976), 280--282.

\bibitem[LM10]{Lawvere-Menni}
F.~William Lawvere and Mat{\'\i}as Menni, \emph{The {H}opf algebra of
  {M}{\"o}bius intervals}, Theory Appl. Categ. \textbf{24} (2010), 221--265.

\bibitem[LMvW13]{Luoto-Mykytiuk-vanWilligenburg}
Kurt Luoto, Stefan Mykytiuk, and Stephanie van Willigenburg, \emph{An
  introduction to quasisymmetric {S}chur functions: {H}opf algebras,
  quasisymmetric functions, and {Y}oung composition tableaux}, SpringerBriefs
  in Mathematics, Springer New York, 2013. \MR{3097867}

\bibitem[Lur09]{Lurie:HTT}
Jacob Lurie, \emph{Higher topos theory}, Annals of Mathematics Studies, vol.
  170, Princeton University Press, Princeton, NJ, 2009. \MR{2522659}

\bibitem[Lur26]{kerodon}
Jacob Lurie, \emph{Kerodon}, \url{https://kerodon.net}, 2026.

\bibitem[Mas19]{Mason:1810.07148}
Sarah~K. Mason, \emph{Recent trends in quasisymmetric functions}, Recent Trends
  in Algebraic Combinatorics (H.~Barcelo, G.~Karaali, and R.~Orellana, eds.),
  Springer International Publishing, Cham, 2019, pp.~239--279.

\bibitem[MN10]{Mastnak-Nica:TAMS2010}
Mitja Mastnak and Alexandru Nica, \emph{Hopf algebras and the logarithm of the
  {$S$}-transform in free probability}, Trans. Amer. Math. Soc. \textbf{362}
  (2010), no.~7, 3705--3743. \MR{2601606}

\bibitem[MR95]{Malvenuto-Reutenauer}
Claudia Malvenuto and Christophe Reutenauer, \emph{Duality between
  quasi-symmetric functions and the {S}olomon descent algebra}, J. Algebra
  \textbf{177} (1995), no.~3, 967--982. \MR{1358493}

\bibitem[OR20]{OzornovaRovelli:MSINCPSSSPSS}
Viktoriya Ozornova and Martina Rovelli, \emph{Model structures for
  {$(\infty,n)$}-categories on (pre)stratified simplicial sets and
  prestratified simplicial spaces}, Algebr. Geom. Topol. \textbf{20} (2020),
  no.~3, 1543--1600. \MR{4105558}

\bibitem[Rot64]{Rota:Moebius}
Gian-Carlo Rota, \emph{On the foundations of combinatorial theory. {I}.
  {T}heory of {M}\"obius functions}, Z. Wahrscheinlichkeitstheorie und Verw.
  Gebiete \textbf{2} (1964), 340--368. \MR{174487}

\bibitem[RV22]{RiehlVerity:EICT}
Emily Riehl and Dominic Verity, \emph{Elements of {$\infty$}-category theory},
  Cambridge Studies in Advanced Mathematics, vol. 194, Cambridge University
  Press, Cambridge, 2022. \MR{4354541}

\bibitem[Sta72]{Stanley:1972}
Richard Stanley, \emph{Ordered structures and partitions}, Memoirs Amer. Math.
  Soc. \textbf{119} (1972).

\bibitem[Sta99]{Stanley:volII}
Richard~P. Stanley, \emph{Enumerative combinatorics. {V}ol. 2}, Cambridge
  Studies in Advanced Mathematics, vol.~62, Cambridge University Press,
  Cambridge, 1999.

\bibitem[Ste21]{Stern:1905.06671}
Walker~H. Stern, \emph{$2$-{S}egal objects and algebras in spans}, J. Homotopy
  Rel Struct. \textbf{16} (2021), 297--361.

\bibitem[SW16]{Shareshian-Wachs:1405.4629}
John Shareshian and Michelle~L. Wachs, \emph{Chromatic quasisymmetric
  functions}, Adv. Math. \textbf{295} (2016), 497--551.

\bibitem[Swe69]{Sweedler}
Moss~E. Sweedler, \emph{Hopf algebras}, W.A. Benjamin, Inc., New York, 1969.

\bibitem[SY19]{SchlankYanovski:ICEHA}
Tomer~M. Schlank and Lior Yanovski, \emph{The {$\infty$}-categorical
  {E}ckmann-{H}ilton argument}, Algebr. Geom. Topol. \textbf{19} (2019), no.~6,
  3119--3170. \MR{4023337}

\bibitem[Ver08]{Verity:WCS1BHT}
Dominic R.~B. Verity, \emph{Weak complicial sets. {I}. {B}asic homotopy
  theory}, Adv. Math. \textbf{219} (2008), no.~4, 1081--1149. \MR{2450607}

\end{thebibliography}

\end{document}